\documentclass{amsart}
\usepackage{amssymb,amsmath,amsthm,epic,eepic,ecltree,pb-diagram,array}

\theoremstyle{plain}
\newtheorem{theorem}{Theorem}[section]
\newtheorem{lemma}[theorem]{Lemma}
\newtheorem{proposition}[theorem]{Proposition}
\newtheorem{definition}[theorem]{Definition}
\newtheorem{corollary}[theorem]{Corollary}

\theoremstyle{definition}

\newtheorem*{scholium}{Scholium}

\newcolumntype{C}{>{$}c<{$}}

\font\wncyrt=wncyr10\newcommand\shadow{\text{\wncyrt sh}} 

\newcommand{\bundleunder}[1]{
  \drawwith{\drawwith{\drawwith{\drawwith\drawline\dottedline{3}}\drawline}\drawline}
  \begin{bundle}{$#1$}
  \chunk{$#11$}
  \chunk{$#12$}
  \chunk{$\cdots$}
  \chunk{$#1q$}
  \end{bundle}}
\newcommand\mathem{\textbf}

\newcommand\N{{\mathbb N}}
\newcommand\R{{\mathbb R}}
\newcommand\Z{{\mathbb Z}}
\newcommand\tree{{\mathcal T}}
\newcommand\treeq{{\tree^{(q)}}}
\newcommand\Omegah{{\widehat\Omega}}
\newcommand\Omegar{{\Omega^{(r)}}}
\newcommand\Omegabr{{\Omega^{[r]}}}
\newcommand\Stab{\mathsf{Stab}}
\newcommand\Aut{\mathsf{Aut}}
\newcommand\Epi{\mathsf{Epi}}
\newcommand\Ker{\mathsf{Ker}}

\begin{document}
\vspace*{1in}
\title{On the Word and Period Growth of some\\ Groups of Tree Automorphisms}
\author{Laurent Bartholdi \and Zoran \v{S}uni\'k}

\keywords{Groups, Rooted Trees, Periodic Groups, Intermediate Growth}
\thanks{The first author expresses his thanks to the ``Swiss National
    Science Foundation''}

\subjclass{\parbox[t]{0.5\textwidth}{%
    \textbf{20F50} (Periodic groups; locally finite groups),\\
    \textbf{16P90} (Growth rate),\\
    \textbf{20E08} (Groups acting on trees)}}
\maketitle 
\thispagestyle{empty} 

    \bigskip 
    \parbox{.48\linewidth}{Section de Math\'ematiques\\
    Universit\'e de Gen\`eve, CP 240\\ 
    1211 Gen\`eve 24, Switzerland\\
    \texttt{laurent.bartholdi@math.unige.ch}}
    \parbox{.48\linewidth}{Department of Mathematical Sciences\\
    SUNY Binghamton\\
    Binghamton, NY 13902, USA\\
    \texttt{shunka@math.binghamton.edu}} \bigskip

\begin{abstract}
  We generalize a class of groups defined by Rostislav Grigorchuk
  in~\cite{grigorchuk:gdegree} to a much larger class of groups,
  and provide upper and lower bounds for their word growth (they are
  all of intermediate growth) and period growth (under a small
  additional condition, they are periodic).
\end{abstract}

\renewcommand\thepart{\Roman{part}}
\renewcommand\labelenumi{(\theenumi)}

%----------------------------------------------------------------
\section{Introduction}
Since William Burnside's original question (``do there exist finitely
generated infinite groups all of whose elements have finite order?''),
dozens of problems have received sometimes unexpected light from the
theory of groups acting on rooted trees: to name a few, John Milnor's
famous Problem~5603~\cite{milnor:5603} on growth, the theory of
just-infinite groups~\cite{grigorchuk:jibg}, of groups of bounded
width~\cite{bartholdi-g:lie} where a sporadic type of group with an
unusual Lie algebra was discovered, the existence of a finitely
presented amenable but non-elementary amenable
group~\cite{grigorchuk:bath}, of groups whose spectrum is a Cantor
set~\cite{bartholdi-g:spectrum}, etc.

The main examples fall roughly in two classes, the \mathem{Grigorchuk
  groups} which have good combinatorial properties (like an ubiquitous
``shortening lemma''), and the \mathem{GGS groups} (named after Slava
Grigorchuk, Narain Gupta and Said Sidki), which have a richer
group-theoretical potential (see for instance~\cite{bartholdi-g:spectrum} 
where torsion-free and torsion groups cohabit). We propose an extension of 
the class defined by R.\ Grigorchuk, and initiate a systematic approach of 
these new groups, which we propose to call \mathem{spinal groups}, since 
the generators are tree automorphisms that are trivial except in the 
neighborhood of a ``spine''. Our hope is that this class is 
\begin{itemize}
\item large enough so that it remains a trove of new examples for
  yet-to-conceive questions, and
\item small enough so that it remains amenable to quantitative
  analysis, in particular thanks to a ``shortening lemma'' that allows
  simple inductive proofs.
\end{itemize}

This paper is roughly comprised of two parts. The first describes
these groups and the combinatorial tools required to fathom them. The
second expands on estimates of the word growth and period growth for
these groups.

The main tools used in the analysis of spinal groups are:
\begin{itemize}
\item A ``shortening lemma''. Each element $g$ in a spinal group
  can be expressed as $g=(g_1,\dots,g_r)h$, where $h$ belongs to a
  finite group, and each of the $g_i$ belongs to a (possibly
  different) spinal group. A lemma (Lemma~\ref{<eta}) states that there 
  is a norm on each group such that the sum of the norms of the $g_i$'s is 
  substantially less than the norm of $g$.
\item A ``portrait representation''. Each element $g$ in a spinal
  group can be described by a subtree $\iota(g)$ of the tree on which
  the group acts, with decorations on nodes of the subtree. For the
  flavor of portraits we use, the subtree $\iota(g)$ is finite, and
  its depth, size etc. carry valuable information on $g$.
\end{itemize}

Each spinal group $G_\omega$ is defined by an infinite sequence
$\omega=\omega_1\omega_2\dots$ of group epimorphisms between two fixed
finite groups. When various conditions are imposed on $\omega$, it is
possible to give good bounds on the growth functions. 

The spinal groups introduced in this paper differ from the Grigorchuk 
examples (see \cite{grigorchuk:burnside},  
\cite{grigorchuk:gdegree} and~\cite{grigorchuk:pgps} for the 
Grigorchuk examples or see the description below in Subsection 
\ref{grigorchuk_examples}) in that they are groups of tree automorphisms 
where the degree is arbitrary (not a prime as in the Grigorchuk examples) 
and the root part of the group does not have to be cyclic (it does not 
have to be abelian either). 

Just to mention a few results that are obtained: 
\begin{itemize}
\item All spinal groups have intermediate growth, more precisely
  subexponential growth, and growth at least $e^{\sqrt n}$. Various
  upper bounds are provided in case the defining sequence $\omega$
  shows some signs of cooperation. We consider two situations:
  $r$-homogeneous and $r$-factorable sequences (see
  Section~\ref{sec:homo}), giving different bounds of the form
  $e^{n^\beta}$.
\item All spinal groups defined through a regular root action are
  periodic. As above, upper and lower bounds for the period growth
  function are provided in the favorable cases; they are both
  polynomial (see Section~\ref{sec:period}).
\item For every $\beta\in(1/2,1)$ there exist spinal groups (even
  among the Grigorchuk examples) whose degree of growth $\gamma_G$ is
  between $e^{n^{\beta}}$ and $e^n$, i.e.\ the degree satisfies 
  $e^{n^{\beta}} \precnsim \gamma_G \precnsim e^n$ (see 
  Theorem~\ref{thm:wdlowerbd}).
\item There exist spinal groups (even among the Grigorchuk examples)
  with at least linear degree of period growth 
  (see Theorem~\ref{thm:periodlb}).
\item The degree of period growth of the first Grigorchuk group is at
  most $n^{3/2}$ (see Theorem~\ref{thm:period ub for 2-groups}).
\end{itemize}

%****************************************************************
\part{Spinal Groups}
\section{Weight Functions, Word and Period Growth}
Let $S=\{s_1,\dots,s_k\}$ be a non-empty set of symbols. A \mathem{weight 
function} on $S$ is any function $\tau:S\rightarrow \R_{>0}$ (note that 
the values are strictly positive).  The \mathem{weight} of any word over 
$S$ is then defined by the extension of $\tau$ to a function, still 
written $\tau:S^*\to\R_{\ge0}$, on the free monoid $S^*$ of words over $S$ 
(note that the empty word is the only word mapped to 0). Let $G$ be an 
infinite group and $\rho:S^*\twoheadrightarrow G$ a surjective monoid 
homomorphism. (Equivalently, $G$ is finitely generated and 
$\rho(S)=\{\rho(s_1),\dots,\rho(s_k)\}$ generates $G$ as a monoid.)  The 
\mathem{weight} of an element $g$ in $G$ with respect to the triple 
$(S,\tau,\rho)$ is, by definition, the smallest weight of a word $u$ in 
$S^*$ that represents $g$, i.e.\ the smallest weight of a word in 
$\rho^{-1}(g)$. The weight of $g$ with respect to $(S,\tau,\rho)$ is 
denoted by $\partial_G^{(S,\tau,\rho)}(g)$. 

For $n$ non-negative real number, the elements in $G$ that have weight at 
most $n$ with respect to $(S,\tau,\rho)$ constitute the \mathem{ball} of 
radius $n$ in $G$ with respect to $(S,\tau,\rho)$, denoted by 
$B_G^{(S,\tau,\rho)}(n)$. The number of elements in 
$B_G^{(S,\tau,\rho)}(n)$ is finite and is denoted by 
$\gamma_G^{(S,\tau,\rho)}(n)$. The function $\gamma_G^{(S,\tau,\rho)}$, 
defined on the non-negative real numbers, is called the \mathem{word 
growth} (or just \mathem{growth}) function of $G$ with respect to 
$(S,\tau,\rho)$. 

If, in addition, $G$ is a torsion group, the following definitions also 
make sense. For $n$ non-negative real number, the maximal order of an 
element in the ball $B_G^{(S,\tau,\rho)}(n)$ is finite and will be denoted 
by $\pi_G^{(S,\tau,\rho)}(n)$. The function $\pi_G^{(S,\tau,\rho)}$, 
defined on the non-negative real numbers, is called the \mathem{period 
growth} function of $G$ with respect to $(S,\tau,\rho)$. 

A partial order $\precsim$ is defined on the set of non-decreasing
functions on $\R_{\geq 0}$ by $f \precsim g$ if there exists a positive
constant $C$ such that $f(n)\leq g(Cn)$ for all $n\in\R_{\geq 0}$.  An
equivalence relation $\sim$ is defined by $f \sim g$ if $f \precsim g$ and 
$g \precsim f$. The equivalence class of $\gamma_G^{(S,\tau,\rho)}$ is 
called the \mathem{degree of growth} of $G$ and it does not depend on the 
(finite) set $S$, the weight function $\tau$ defined on $S$ and the 
homomorphism $\rho$. The equivalence class of $\pi_G^{(S,\tau)}$ is called 
the \mathem{degree of period growth} of $G$ and it also does not depend on 
the triple $(S,\tau,\rho)$. 

Of course, when we define a weight function on a group $G$ we usually pick 
a finite generating subset of $G$ closed for inversion and not containing 
the identity, assign a weight function to those generating elements and 
extend the weight function to the whole group $G$ in a natural way, thus 
blurring the distinction between a word over the generating set and the 
element in $G$ represented by that word and completely avoiding the 
discussion of $\rho$. In most cases everything is still clear that way. 

Let us mention that the standard way to assign a weight function is to 
assign the weight 1 to each generator. In that case we denote the weight 
of a word $u$ by $|u|$ and call it the \mathem{length} of $u$. In this 
setting, the length of the group element $g$ is the distance from $g$ to 
the identity in the Cayley graph of the group.  

Since the degree of growth is invariant of the group we are more 
interested in it than in the actual growth function for a given generating 
set. For any finitely generated infinite group $G$, the following 
trichotomy exists: $G$ is of 
\begin{itemize}
\item \mathem{polynomial growth} if $\gamma_G(n)\precsim n^d$ for some
  $d\in\N$;
\item \mathem{intermediate growth} if
  $n^d\precnsim\gamma_G(n)\precnsim e^n$ for all $d\in\N$;
\item \mathem{exponential growth} if $e^n\sim\gamma_G(n)$.
\end{itemize}
We also say $G$ is of \mathem{subexponential growth} if
$\gamma_G(n)\precnsim e^n$ and of \mathem{superpolynomial growth} if
$n^d\precnsim\gamma_G(n)$ for all $d\in\N$.

The classes of groups of polynomial and exponential growth are clearly
non-empty: the former consists, by a theorem of Mikhail
Gromov~\cite{gromov:nilpotent}, precisely of virtually nilpotent
groups and the latter contains, for instance, all non-elementary
hyperbolic groups~\cite{ghys-h:gromov}.

As a consequence of Gromov's theorem, if a group $G$ is of polynomial
growth then its growth function $\gamma_G$ is equivalent to $n^d$ for
an integer $d$. There even is a formula giving $d$ in terms of the
lower central series of $G$, due to Yves Guivarc'h and Hyman
Bass~\cite{guivarch:poly1,bass:nilpotent}.

By Tits' alternative~\cite{tits:linear}, there are no examples of groups 
of intermediate growth among the linear groups. However, R.\ Grigorchuk 
discovered examples of groups of intermediate growth by studying piecewise 
diffeomorphisms of the real line~\cite{grigorchuk:growth}, and other 
examples followed~\cite{gupta-f:growth2,bartholdi:ggs}. 

%-----------------------------------------------------------------------------
\section{The Groups}
The class of groups we are about to define is a generalization of the
class of Grigorchuk $p$-groups introduced
in~\cite{grigorchuk:burnside},~\cite{grigorchuk:gdegree}
and~\cite{grigorchuk:pgps}. An intermediate generalization was already
suggested by Grigorchuk in~\cite{grigorchuk:pgps}, but seems never to
have been pursued.

Also, Alexander Rozhkov gives even more general constructions of
similar type in~\cite{rozhkov:aleshin}.

In the original description, Grigorchuk groups are given as groups of
permutations of the unit interval from which a set of measure 0 is
removed. In this paper we find it more convenient to describe the
groups as groups of automorphisms of the $q$-regular rooted tree.

\subsection{Infinite regular rooted trees and tree automorphisms}
The approach we take here follows~\cite{brin:igt} and~\cite{harpe:cgt}.

Fix once and for all an integer $q\ge2$, and set $Y=\{1,2,\dots,q\}$.
We think of the \mathem{$q$-regular rooted tree} $\treeq$ as the set
of finite words over $Y$ related by the prefix ordering. Recall that a
word over the alphabet $Y$ is just a finite sequence of elements of
$Y$; for convenience we start the indexing of the letters at $0$. In
the prefix ordering, $u \leq v$ if and only if $u$ is a prefix of $v$,
and an edge joins two vertices in $\treeq$ precisely when one is an
immediate successor of the other.

Every finite word represents a vertex of the tree: the empty word
represents the \mathem{root}, the words $1,2,\dots,q$ represent the
vertices on the first level below the root, the two-letter words
$11,12,\dots,1q$ represent the vertices on the second level below the
vertex $1$, etc.

\vspace{5mm} \centerline{ 
\drawwith{\drawwith{\drawwith{\drawwith\drawline\dottedline{3}}\drawline}\drawline} 
\begin{bundle}{$\emptyset$}
\chunk{\bundleunder{1}} \chunk{\bundleunder{2}} \chunk{$\cdots$} 
\chunk{\bundleunder{q}} 
\end{bundle}}
\vspace{3mm} \centerline{The tree $\treeq$} \vspace{5mm}

The vertex $u$ is above the vertex $v$ in the tree if and only if $u$
is a prefix of $v$. The vertex $v$ is a child of the vertex $u$ if and
only if $v=ui$ for some $i\in Y$. The words of length $k$ constitute
the \mathem{level} $L_k$ in the tree.

An automorphism of the tree $\treeq$ is any permutation of the
vertices in $\treeq$ preserving the prefix ordering (and therefore,
also the length). Every automorphism $g$ of $\treeq$ induces a
permutation of the set $\partial\treeq=Y^\N$ of infinite sequences
(again indexed from $0$) over $Y$ in a natural way. Geometrically,
$\partial\treeq$ is the boundary of $\treeq$. For two infinite
sequences $u$ and $v$ in $\partial\treeq$ we define $u \wedge v$ to
be the longest common prefix of $u$ and $v$. An automorphism $g$ of
the tree $\treeq$ induces a permutation $\bar{g}$ of $\partial\treeq$
satisfying
\begin{equation} \label{prefcond}
 |\bar{g}(u) \wedge \bar{g}(v)| = |u \wedge v|
\end{equation}
for all infinite sequences $u$, $v$ in $\partial\treeq$. Conversely,
every permutation $\bar{g}$ of $\partial\treeq$
satisfying~(\ref{prefcond}) induces an automorphism $g$ of the tree
$\treeq$ in a natural way.

In the sequel, it will be convenient for us to define some tree
automorphisms by using this alternative way (permutations of infinite
sequences). Actually, we will not distinguish between the two ways at
all and we will switch back and forth between the two points of view.
Also, from now on we will write $\tree$ instead of $\treeq$, and will
denote the automorphism group of $\tree$ by $\Aut(\tree)$.

For a word $u$ over $Y$ denote by $\tree_u$ the set of words in
$\tree$ that have $u$ as a prefix. The set $\tree_u$ has a tree
structure for the prefix ordering and it is isomorphic to $\tree$ by
the canonical isomorphism deleting the prefix $u$. Any automorphism
$g$ of $\tree$ that fixes the word $u$ induces an automorphism
$g_{|u}$ of $\tree_u$ by restriction. Every automorphism $g_{|u}$ of
$\tree_u$ fixes $u$ and induces an automorphism $g_u$ of $\tree$,
which acts on the word $w$ exactly as $g_{|u} \in \Aut(\tree_u)$ acts
on the $w$ part of the word $uw$, namely by $ug_u(w)=g_{|u}(uw)$.  The
map $\varphi_u$ defined by $g\mapsto g_{|u}\mapsto g_u$ is a
surjective homomorphism from the stabilizer $\Stab(u)$ of $u$ in
$\Aut(\tree)$ to the automorphism group $\Aut(\tree)$.

Let $\Stab(L_1)$ be the stabilizer of the first level of $\tree$ in
$\Aut(\tree)$, i.e.\ $\Stab(L_1)=\bigcap_{i=1}^q \Stab(i)$. The homomorphism 
$\psi:\Stab(L_1) \rightarrow \Pi_{i=1}^q \Aut(\tree)$ given by 
\[ \psi(g)=(\varphi_1(g),\varphi_2(g),\dots,\varphi_q(g))=(g_1,g_2,\dots,g_q)\]
is an isomorphism.

Similarly, let $\Stab(L_r)$ be the stabilizer of the $r$-th level of
$\tree$ in $\Aut(\tree)$:
\[\Stab(L_r)=\bigcap\{\Stab(u)|\,u\text{ is a $r$-letter word
  in }\tree\}.\] 
The homomorphism $\psi_r:\Stab(L_r) \rightarrow\prod_{i=1}^{q^r} \Aut(\tree)$ 
given by 
\[ \psi_r(g)=(\varphi_{1\dots11}(g),\varphi_{1\dots12}(g),\dots,\varphi_{q\dots qq}(g))= 
    (g_{1\dots 11},g_{1\dots 12},\dots,g_{q\dots qq}) \]
is an isomorphism.

%-----------------------------------------
\subsection{The construction of the groups}
Let $G_A$ be a group (called the \mathem{root group}) acting faithfully 
and transitively on $Y$ (therefore, $G_A$ is finite of order at least $q$ 
and most $q!$). Further, let $G_B$ be a finite group (called the 
\mathem{level group}) such that the set $\Epi(G_B, G_A)$ of surjective 
homomorphisms from $G_B$ to $G_A$ is non-empty. When an epimorphism in 
$\Epi(G_B,G_A)$ is called $\omega_i$, denote the kernel $\Ker(\omega_i)$ by 
$K_i$. We impose additional requirements on $G_B$ by asking that 
\begin{itemize}
\item the union of all these kernels is $G_B$ (so that every element
  in $G_B$ is sent to the identity by some homomorphism in
  $\Epi(G_B,G_A)$);
\item their intersection is trivial (which, among the other things,
  says that $G_B$ is a subdirect product of several copies of $G_A$).
\end{itemize}

The set $\Omegah$ is defined as the set of infinite sequences $\omega = 
\omega_1\omega_2\dots$ over $\Epi(G_B,G_A)$ such that every non-trivial 
element $g$ of $G_B$ both appears and does not appear in infinitely many 
of the kernels $K_1,K_2,\dots$. Note that the indexing of the sequences in 
$\Omegah$ starts with $1$. Equivalently, we might say that $\Omegah$ 
consists of the sequences $\omega = \omega_1\omega_2\dots$ over 
$\Epi(G_B,G_A)$ such that for every $i$ we have 
\[ \bigcup_{i\leq j} K_j = G_B \qquad \text{and} \qquad  
\bigcap_{i\leq j} K_j = 1. \] 
It is true that $\Omegah$ depends on $G_B$ and $G_A$, but we will avoid 
any notation emphasizing that fact. The \mathem{shift operator} 
$\sigma:\Omegah \rightarrow \Omegah$ is defined by 
$\sigma(\omega_1\omega_2\dots) = \omega_2\omega_3\dots$. 

The root group $G_A$ acts faithfully on the boundary of $\tree$ by acting 
on the $0$-coordinate in $Y^\N=\partial\tree$, namely by 
\[ g(y_0y_1y_2\dots) = g(y_0)y_1y_2\dots.\]
The automorphism of $\tree$ induced by $g \in G_A$ will also be denoted by 
$g$ and the set of non-identity automorphisms of $\tree$ induced by $G_A$ 
will be denoted by $A$. Letters like $a,a_1,a',\dots$ are reserved for the 
elements in $A$.   

Given a sequence $\omega$ in $\Omegah$ we define an action of the level 
group $G_B$ on the tree $\tree$ as follows: 
\[ \begin{align*}
  g(q\dots q1y_{n+1}y_{n+2}\dots) &=
  q\dots q1\omega_{n+1}(g)(y_{n+1})y_{n+2}\dots;\\
  g(y) &= y\text{ for any word $y$ not starting with }q\dots q1.
\end{align*} \]
The group $G_B$ acts faithfully on $\tree$ as a group of tree 
automorphisms.  The tree automorphism corresponding to the action of $g$ 
will be denoted by $g_\omega$ and the set of non-identity tree 
automorphisms induced by $G_B$ will be denoted by $B_\omega$. The abstract 
group $G_B$ is canonically isomorphic to the group of tree automorphisms 
$G_{B_\omega}$ for any $\omega$ so that we will no make too much 
difference between them and will frequently omit the index $\omega$ in the 
notation. The index will be omitted in $B_\omega$ as well. Letters like 
$b,b_1,b',\dots$ are reserved for the elements in $B$.   

What really happens is that $a \in G_A$ acts at the root of $\tree$ by
permuting the subtrees $\tree_1,\dots,\tree_q$. On the other hand, the
action of $b \in G_B$ is prescribed by $\omega$. Namely, $b$ acts on the 
subtree $\tree_1$ exactly as $\omega_1(b) \in G_A$ would act on $\tree$, 
on the subtree $\tree_{q1}$ exactly as $\omega_2(b) \in G_A$ would act on 
$\tree$, $\dots$; and $b$ acts trivially on subtrees not of the form 
$\tree_{q\dots q 1w}$. 

\vspace{5mm} \centerline{ 
\drawwith{\drawwith{\drawwith{\drawwith\drawline\dottedline{3}}\drawline}\drawline} 
\begin{bundle}{}
  \chunk{$\omega_1(b)$} \chunk{$1$} \chunk{$\cdots$}
  \chunk{\begin{bundle}{}
      \chunk{$\omega_2(b)$} \chunk{$1$} \chunk{$\cdots$}
      \chunk{\begin{bundle}{} \chunk{$\omega_3(b)$} \chunk{$1$} 
            \chunk{$\dots$} \chunk{} 
            \end{bundle}}
    \end{bundle}}
\end{bundle}}
\vspace{3mm} \centerline{The automorphism $b_\omega \in G_B$} \vspace{5mm} 

We define now the main object of our study:

\begin{definition} 
  For any sequence $\omega \in \Omegah$, the subgroup of the
  automorphism group $\Aut(\tree)$ of the tree $\tree$ generated by $A$
  and $B_\omega$ is denoted by $G_\omega$ and called the
  \mathem{spinal group} defined by the sequence $\omega$.
\end{definition}

Note that we could define a still larger class of groups if we avoided
some of the self-imposed restrictions above. However, many of the
properties that follow and that interest us would not hold in that
larger class.

The groups $G_\omega$, like all groups acting on a rooted tree, can be
described as \mathem{automata
  groups}~\cite{gecseg-c:ata,bartholdi-g:spectrum}, where the elements
of the group are represented by Mealy (or Moore) automata with
composition of automata as group law. These automata will be finite
automata if and only if the sequence $\omega$ is regular, i.e.\
describable by a finite automaton. We will not pursue this topic here.

%-----------------------------------------
\subsection{The examples of Grigorchuk} \label{grigorchuk_examples}
Before we proceed with our investigation of the constructed groups let us 
describe exactly what groups were introduced by Grigorchuk 
in~\cite{grigorchuk:burnside},~\cite{grigorchuk:gdegree} 
and~\cite{grigorchuk:pgps}. Each Grigorchuk $p$-group acts on a rooted 
$p$-regular tree, for $p$ a prime, the root group $G_A \cong \Z/p\Z$ is 
the group of cyclic permutations of $Y=\Z/p\Z=\{1,2,\dots,p\}$ generated 
by the cyclic permutation $a=(1,2,\dots,p)$, the level group $G_B$ is 
isomorphic to $\Z/p\Z \times\Z/p\Z$ and only the following $p+1$ 
homomorphisms from $G_B$ to $G_A$ are used in the construction of the 
infinite sequences in $\Omegah$. These epimorphisms are written 
$\left[\begin{smallmatrix}u\\ v\end{smallmatrix}\right]$ to mean the 
linear functionals on $\Z/p\Z\times\Z/p\Z$ given by $(x,y)\mapsto ux+vy$: 

\[ \begin{bmatrix} 1 \\ 0 \end{bmatrix},
 \begin{bmatrix}1 \\ 1 \end{bmatrix}, 
\begin{bmatrix} 1 \\ 2 \end{bmatrix}, \dots,
\begin{bmatrix} 1 \\ p-1 \end{bmatrix},
\begin{bmatrix} 0 \\ 1 \end{bmatrix}. \]

The most known and investigated example is the \mathem{first Grigorchuk 
group}~\cite{grigorchuk:burnside}, which is defined as above for $p=2$ 
where the sequence 
\[ \omega=\begin{bmatrix} 1 \\ 1 \end{bmatrix}\begin{bmatrix} 1 \\ 0 
\end{bmatrix}\begin{bmatrix} 0 \\ 1 \end{bmatrix}\begin{bmatrix} 1 \\ 1
\end{bmatrix}\begin{bmatrix} 1 \\ 0 \end{bmatrix}\begin{bmatrix} 0 \\ 1
\end{bmatrix}\dots \]
is periodic of period $3$. 

In the case of Grigorchuk $2$-groups (acting on the binary tree), it is 
customary to denote the only nontrivial element of the root group 
$G_A=\Z/2\Z$ by $a$ and the three nontrivial elements of the level group 
$G_B=\Z/2\Z \times \Z/2\Z$ by $b$, $c$ and $d$.  There are only three 
epimorphisms from $G_B$ to $G_A$ and each of them maps exactly one of the 
$B$-generators $b$, $c$, $d$ to $1$ and the other two to $a$. The 
epimorphisms sending $d$, $c$ and $b$, respectively, to $1$ are denoted by 
$0$, $1$ and $2$. Then the set of admissible sequences $\Omegah$ consists 
of all those sequences that contain each of these three epimorphisms 
infinitely many times, i.e.\ sequences over $\{0,1,2\}$ that have 
infinitely many appearances of each of the letters $0$, $1$ and~$2$. In 
this terminology, the first Grigorchuk group is defined by the sequence 
$012012012\dots$. 

%-----------------------------------------
\subsection{More examples} 
It is not difficult to construct many examples where $G_A$ and $G_B$ are 
abelian, but the construction of examples where $G_A$ and $G_B$ are not 
abelian is not obvious. The following example, which allows different 
generalizations, was suggested by Derek Holt. 

Let $G_B= \langle b_1,b_2,b_3,b_4,b_5,b_6,x_{12},x_{34}\rangle$ where 
$b_1,b_2,b_3,b_4,b_5,b_6$ all have order 3 and commute with each other, 
$x_{12}$ and $x_{34}$ have order 2 and commute and 
\[ b_i^{x_{jk}} = 
    \begin{cases} b_i, \qquad \text{if} \; i\in\{j,k\}\\
                  b_i^{-1}, \qquad \text{otherwise} 
    \end{cases}. \]
In other words $G_B$ is the semidirect product $(\Z/3\Z)^6 \rtimes 
(\Z/2\Z)^2$ where $(\Z/2\Z)^2=\langle x_{12},x_{34} \rangle$ and $x_{12}$ 
fixes the first two coordinates of $(\Z/3\Z)^6$ and acts by inversion on 
the last $4$, $x_{34}$ fixes the middle $2$ coordinates and acts by inversion 
on the other $4$ and, consequently, $x_{56}=x_{12}x_{34}$ fixes the last two 
coordinates and inverts the first $4$. 

The following $12$ subgroups are normal in $G_B$, their intersection is 
trivial, their union is $G_B$, and each factor is isomorphic to the 
symmetric group $\Z/3Z \rtimes \Z/2\Z = S_3$ which is then taken to be 
$G_A$: 
\[\begin{matrix}
\langle b_1,b_3,b_4,b_5,b_6,x_{12} \rangle,&
\langle b_1,b_2,b_3,b_5,b_6,x_{34} \rangle,&
\langle b_1,b_2,b_3,b_4,b_5,x_{56} \rangle, \\
\langle b_2,b_3,b_4,b_5,b_6,x_{12} \rangle,&
\langle b_1,b_2,b_4,b_5,b_6,x_{34} \rangle,&
\langle b_1,b_2,b_3,b_4,b_6,x_{56} \rangle, \\
\langle b_1b_2,b_3,b_4,b_5,b_6,x_{12} \rangle,&
\langle b_1,b_2,b_3b_4,b_5,b_6,x_{34} \rangle,&
\langle b_1,b_2,b_3,b_4,b_5b_6,x_{56} \rangle, \\
\langle b_1b_2^2,b_3,b_4,b_5,b_6,x_{12} \rangle,&
\langle b_1,b_2,b_3b_4^2,b_5,b_6,x_{34} \rangle,&
\langle b_1,b_2,b_3,b_4,b_5b_6^2,x_{56} \rangle.
\end{matrix}\]
We consider the $12$ epimorphisms from $G_B$ to $G_A$ that are the
quotient maps by these normal subgroups, and accept in $\Omegah$ all
sequences $\omega$ that uses each of these $12$ homomorphisms infinitely
often.

%-----------------------------------------------------------------------------
\section{Some Tools for Investigation of the Groups} 
In this section we introduce the tools and constructions we will use
in the investigation of the groups along with some basic properties
that follow quickly from the given considerations.

\subsection{Triangular weights and minimal forms}
The finite set $S_\omega = A \cup B_\omega$ is the \mathem{canonical
  generating set} of $G_\omega$. The generators in $A$ are called
\mathem{$A$-generators} and the generators in $B_\omega$ are called
\mathem{$B$-generators}. Note that $S_\omega$ does not contain the
identity and generates $G_\omega$ as a monoid, since it is closed
under inversion.

A weight function $\tau$ on $S$ will be called \mathem{triangular} if
\[ \tau(a_1) + \tau(a_2) \geq \tau(a_1a_2) \qquad \text{and} \qquad
    \tau(b_1) + \tau(b_2) \geq \tau(b_1b_2), \] 
for all $a_1,a_2\in A$ and $b_1,b_2 \in B$ such that $a_1a_2\in A$ and 
$b_1b_2\in B$.

Every $g$ in $G_\omega$ admits a minimal form with respect to a
triangular weight $\tau$
\begin{equation} \label{form}
[a_0]b_1a_1b_2a_2 \dots a_{k-1}b_k[a_k]
\end{equation}
where all $a_i$ are in $A$ and all $b_i$ are in $B$, and $a_0$ and
$a_k$ are optional. This is clear, since the appearance of two
consecutive $A$-letters can be replaced either by the empty word (if
the corresponding product in $G_\omega$ is trivial) or by another
$A$-letter (if the product corresponds to a non-trivial element in
$G_\omega$). In each case the reduction of this type does not increase
the weight, while it decreases the length. The same argument is valid
for consecutive $B$-letters.

Note that the standard weight function, the length, is triangular and,
therefore, admits a minimal form of type~(\ref{form}).

Relations of the following $4$ types:
\[ a_1a_2 \rightarrow 1, \quad a_3a_4 \rightarrow a_5, \quad
b_1b_2 \rightarrow 1, \quad b_3b_4 \rightarrow b_5, \] that follow from 
the corresponding relations in $G_A$ and $G_B$ for $a_i \in A$ and $b_j 
\in B$ are called \mathem{simple relations}. A \mathem{simple 
  reduction} is any single application of a simple relation from left
to right (indicated above by the arrows). Any word of the
form~(\ref{form}) will be called a \mathem{reduced word} and any word
can uniquely be rewritten in reduced form using simple reductions. Of
course, the word and its reduced form represent the same element.

Note that the system of reductions described above is complete, i.e.\
it always terminates with a word in reduced form and the order in
which we apply the reductions does not change the final reduced word
obtained by the reduction.  The second property, known as the
\emph{Church-Rosser} property, is not very important for us since we
can agree to a standard way of performing the reductions (for example,
always reduce at a position as close to the beginning of the word as
possible).

%-------------------------------
\subsection{Some homomorphisms and subgroups}
The intersection $\Stab(L_1)\cap G_\omega$, denoted by $H_\omega$, is a
normal subgroup of $G_\omega$ (since $\Stab(L_1)$ is normal in
$\Aut(\tree)$) and it consists of those elements of $G_\omega$ that fix
the first symbol of each infinite word in $Y^\N=\partial\tree$.

Since each element in $B$ fixes the first level, a word $u$ over $S$
represents an element in $H_\omega$ if and only if the word in
$A$-letters obtained after deleting all the $B$-letters in $u$
represents the identity element.

Further, $H_\omega$ is the normal closure of $B_\omega$ in $G_\omega$,
with $G_\omega/H_\omega \cong G_A$, and $H_\omega$ is generated by the
elements $b_\omega^{g}=gb_\omega g^{-1}$ for $b$ in $B$ and $g$ in $G_A$. 

Denote by $\varphi_i^\omega$ the homomorphism obtained by restricting
$\varphi_i$ to $H_\omega$ in the domain and to the image
$\varphi_i(H_\omega)$ in the codomain and let us calculate this image.
Clearly, $\psi(b_\omega)=(\omega_1(b),1,\dots,1,b_{\sigma\omega})$. For 
any $a$ in $A$, $\psi(b_\omega^a)$ has the same components as 
$\psi(b_\omega)$ does but in different positions depending on $a$. For 
example, if $a$ is the cyclic permutation $(1,2,\dots,q)$ (meaning $1 
\mapsto 2 \mapsto \dots \mapsto q \mapsto 1$), the images of 
$b_\omega^{a^j}$ under various $\varphi_i^\omega$ are given in 
Table~\ref{table:phi}. 
\begin{table}
  \[\begin{tabular}{|C|C C C C C C C C|} \hline
    &\varphi_1^\omega &\varphi_2^\omega &\varphi_3^\omega &\cdots &\varphi_{i+1}^\omega &\cdots &\varphi_{q-1}^\omega &\varphi_q^\omega\\[2pt]\hline
    b_\omega &\omega_1(b) &1 &1 &\cdots &1 &\cdots &1 &b_{\sigma\omega}\\
    b_\omega^{a} &b_{\sigma\omega} &\omega_1(b) &1 &\cdots &1 &\cdots &1 &1\\
    b_\omega^{a^2} &1 &b_{\sigma\omega}& \omega_1(b) &\cdots &1 &\cdots &1 &1\\
    \vdots &\vdots &\vdots &\vdots &\ddots &\vdots &\phantom{\omega_1(b)} &\vdots &\vdots\\
    b_\omega^{a^i} &1 &1 &1 &\cdots &\omega_1(b) &\cdots &1 &1 \\
    \vdots &\vdots &\vdots &\vdots &\phantom{\omega_1(b)} &\vdots &\ddots &\vdots &\vdots\\
    b_\omega^{a^{q-2}} &1 &1 &1 &\cdots &1 &\cdots &\omega_1(b) &1 \\
    b_\omega^{a^{q-1}} &1 &1 &1 &\cdots &1 &\cdots &b_{\sigma\omega} &\omega_1(b)\\ \hline
  \end{tabular}\]
  \vspace{1ex}
  \caption{The maps $\varphi_i$ associated with the permutation
    $a=(1,2,\dots,q)$}\label{table:phi}
\end{table}

Since $\omega_1$ is surjective and the root group acts transitively on
$Y$ we get all $A$ and all $B$-generators in the image of every
$\varphi_i^\omega$.  Therefore $\varphi_i^\omega:H_\omega \rightarrow
G_{\sigma\omega}$ is a surjective homomorphism for all
$i=1,2,\dots,q$.

The homomorphism $\psi^\omega:H_\omega \rightarrow \prod_{i=1}^q
G_{\sigma\omega}$ given by
\[\psi^\omega(g)=(\varphi_1^\omega(g),\dots,\varphi_q^\omega(g)) =
(g_1,g_2,\dots,g_q)\] is a subdirect embedding, i.e.\ is surjective on
each factor. We will avoid the superscript $\omega$ as much as
possible.

Similarly, the intersection $\Stab(L_r)\cap G_\omega$, denoted by
$H_\omega^{(r)}$, is a normal subgroup of $G_\omega$ (since
$\Stab(L_r)$ is normal in $\Aut(\tree)$) and consists of those elements
of $G_\omega$ that fix the first $r$ symbols of each infinite word in
$\partial\tree$. The homomorphism $\psi_r^\omega:H_\omega^{(r)}
\rightarrow \prod_{i=1}^{q^r} G_{\sigma\omega}$ given by
\[ \psi_r^\omega(g)=(\varphi_{1\dots1}(g),\dots,\varphi_{q\dots q}(g)) =  
    (g_{1\dots1},\dots,q_{q\dots q}) \]
is a subdirect embedding.

We end this subsection with a few easy facts, whose proof we omit:
\begin{lemma}\label{easylemma}
  For any $h\in H_\omega$, $g\in G_A$, $b\in B$ and
  $i\in\{1,\dots,q\}$, we have
  \begin{enumerate}
  \item $|h_i|\leq (|h|+1)/2$.
  \item $\varphi_i(h^g) = \varphi_{g^{-1}(i)}(h)$.
  \item The coordinates of $\psi(b^g)$ are: $\omega_1(b)$ at the
    coordinate $g(1)$, $b$ at $g(q)$ and $1$ elsewhere.
  \end{enumerate}
\end{lemma}

\begin{proposition}
  The group $G_\omega$ is infinite for every $\omega$ in $\Omegah$.
\end{proposition}
\begin{proof}
  The proper subgroup $H_\omega$ maps onto $G_{\sigma\omega}$ (for
  instance under $\varphi_1^\omega$), so
  $|G_\omega|\ge|A|\cdot|G_{\sigma\omega}|$. Then the proper subgroup
  $H_{\sigma\omega}$ of $G_{\sigma\omega}$ maps onto
  $G_{\sigma^2\omega}$ (under $\varphi_1^{\sigma\omega}$), so
  $|G_\omega|\ge|A|^2|G_{\sigma^2\omega}|$; etc.
\end{proof}

\begin{proposition}
  The group $G_\omega$ is residually finite.
\end{proposition}
\begin{proof}
  $G_\omega$ is a subgroup of $\Aut(\tree)$, which clearly is
  residually finite: it is approximated by its finite quotients given
  by the action on $Y^n$, for any $n\in\N$.
\end{proof}

\begin{proposition}
  The group $G_\omega$ has a trivial center.
\end{proposition}
\begin{proof}
  First, we will prove that if $g\in G_\omega$ is central then $g$
  must be in $H_\omega$.
  
  Let $g=ha$ where $h\in H_\omega$ and $a\in A$. If $a(1)=1$ and
  $a(i)=j$ for some $i \not= j$ then $g$ does not commute with the
  elements $a'\in A$ such that $a'(1)=i$. If $a(1)\not= 1$ then choose
  $b\in B$ with $b \not \in K_1$ and consider
  \[\varphi_1([g,b])=\varphi_1(gbg^{-1}b^{-1})=
  \varphi_1(h)\varphi_{a^{-1}(1)}(b)\varphi_1(h^{-1})\omega_1(b^{-1}).\]
  It is clear that
  $\varphi_1(h)\varphi_{a^{-1}(1)}(b)\varphi_1(h^{-1}) \in
  H_{\sigma\omega}$ since $a^{-1}(1)\not=1$, which, along with the
  fact that $\omega_1(b^{-1})\not=1$, gives $\varphi_1([g,b]) \not\in
  H_{\sigma\omega}$ and therefore $[g,b]\not=1$.
  
  Now, we proceed by induction on the length of the elements and we
  prove the statement for all $\omega$ simultaneously. From the above
  discussion it is clear that no $A$-generator and no $B$-generator
  outside of $K_1$ is in the center. Consider a generator $b\in K_1$
  and choose an element $a\in A$ with $a(q)=1$. Then
  $\varphi_1([a,b])=\varphi_{a^{-1}(1)}(b)\omega_1(b^{-1})=b\not=1$.
  Therefore no element in $K_1$ is in the center and we have completed
  the basis of the induction.
  
  Consider an element $g\in G_\omega$ of length $\geq 2$. If $g\not\in
  H_\omega$ we already know that $g$ is not in the center. Let $g\in
  H_\omega$. At least one of the projections, say $g_i \in
  G_{\sigma\omega}$, is not trivial. Since $g_i$ has strictly shorter
  length than $g$, we obtain that $g_i$ is not in the center of
  $G_{\sigma\omega}$ so that $g$ is not in the center of $H_\omega$
  (and therefore not in the center of $G_\omega$).
\end{proof}

\begin{proposition}\label{general-finite-subgp}
  The subgroup $D_r = \langle A,K_r\rangle$ of $G_\omega$ is finite
  for any $r\in\N$.
\end{proposition}
\begin{proof}
  Take any $g\in D_r\cap H_\omega^{(r)}$, and consider the coordinates
  of $\psi_r^\omega(g)$. They all belong to $K_r$, whence $D_r\cap
  H_\omega^{(r)}$ is a finite group of order at most $|K_r|^{q^r}$,
  and $D_r$ is finite, since $H_\omega^{(r)}$ is of finite index in
  $G_\omega$.
\end{proof}

\begin{corollary}\label{finite-subgp}
  For $\omega\in\Omegah$, the subgroup $\langle A,b\rangle$ of
  $G_\omega$ is finite for any $b\in G_B$.
\end{corollary}

%--------------------------------
\subsection{Tree decomposition of reduced words} 
The following construction corresponds directly to a construction 
exhibited in~\cite{grigorchuk:gdegree} and~\cite{grigorchuk:pgps}. 
Let 
\[ F = [a_0]b_1a_1b_2a_2 \dots a_{k-1}b_k[a_k] \]
be a reduced word in $S$ representing an element in $H_\omega$. We rewrite the 
element $F$ of $G_\omega$ in the form 
\[ \begin{aligned}
    F &= b_1^{[a_0]}b_2^{[a_0]a_1} \dots b_k^{[a_0]a_1\dots a_{k-1}}
    [a_0]a_1\dots a_{k-1}[a_k] \\
    &= b_1^{g_1}b_2^{g_2} \dots b_k^{g_k}[a_0]a_1\dots a_{k-1}[a_k], 
\end{aligned}\]  
where $g_i=[a_0]a_1\dots a_{i-1} \in G_A$. We know that $[a_0]a_1\dots
a_{k-1}[a_k]=1$, since $F$ is in $H_\omega$. Next, using the definition of 
$\omega_1$ and a table similar to Table~\ref{table:phi} (but for all 
possible $a$) we compute the (not necessarily reduced) words 
$\overline{F_1}, \dots, \overline{F_q}$ representing the elements 
$\varphi_1(F),\dots,\varphi_q(F)$ of $G_{\sigma\omega}$, respectively. 
Then we reduce these $q$ words using simple reductions and obtain the 
reduced words $F_1,\dots,F_q$. We still have $\psi(F)=(F_1,\dots, F_k)$. 
The order in which we perform the reductions is unimportant since the 
system of simple reductions is complete. Thus the rooted $q$-ary labeled 
tree of depth 1 whose root is decorated by the word $F$ and 
its $q$ children by the words $F_1,\dots,F_q$ is well defined and we call 
it the \mathem{depth-$1$ decomposition of $F$}. 

Note that each $B$-letter $b$ from $F$ contributes exactly one
appearance of the letter $b$ to one of the words
$\overline{F_1},\dots,\overline{F_q}$ and, possibly, one $A$-letter to
another word. Thus, the length of any of the reduced words
$F_1,\dots,F_q$ does not exceed $k$ i.e.\ does not exceed $(n+1)/2$
where $n$ is the length of $F$.

Given an $r>1$ and a reduced word $F$ representing an element in 
$H^{(r)}_\omega$, we construct a rooted $q$-ary labeled tree of depth $r$ 
inductively as follows: the root is decorated by $F$ and 
the decompositions of depth $r-1$ of $F_1,\dots,F_q$ are attached to the 
$q$ children of the root. We call this tree the \mathem{depth-$r$ 
  decomposition of $F$}.

Note that the vertices on the second level in the decomposition are
decorated by $F_{11}$, $F_{12}$ $\dots$, $F_{qq}$ where $F_{ij}$ have the 
property 
\[ \psi^{\sigma\omega}(F_i)=
(\varphi_1^{\sigma\omega}(F_i),\varphi_2^{\sigma\omega}(F_i),
\dots,\varphi_q^{\sigma\omega}(F_i))= (F_{i1},F_{i2}\dots,F_{iq}). \] The 
vertices on the third level are decorated by 
$F_{111},F_{112},\dots,F_{qqq}$, etc. 

%--------------------------------------
\subsection{The commutators} 
We now determine the commutator subgroup $[G_\omega, G_\omega]$ of 
$G_\omega$ along with the abelianization $G_\omega^{ab}= 
G_\omega/[G_\omega,G_\omega]$. 

For a word $F$ over $S$, define the word $F_B$ to be the $B$-word
obtained after the removal of all $A$-letters in $F$. Then define the
following set of words over $S$:
\[ \Ker(\rho_B^{ab}) = \{ F\in S^*|\,F_B\text{ represents the identity in }G_B^{ab}\}.\]

\begin{lemma}\label{lemma:welldef'd}
  If $F$ is a word over $S$ representing identity in $G_\omega$, then
  $F \in \Ker(\rho_B^{ab})$. In other words, all relators in $G_\omega$
  come from $\Ker(\rho_B^{ab})$.
\end{lemma}

We may therefore consider the map $p_B^{ab}:G_\omega \rightarrow G_B^{ab}$ 
given by $p_B^{ab}(g)=F_B$, where $F$ is any word over $S$ representing 
$g$, which is clearly surjective, and is well-defined by the lemma above: 
\[ \begin{diagram}
  \node{S^*}\arrow{e,t}{\rho_B^{ab}}\arrow{s,l}{\rho_\omega}\node{G_B^{ab}}\\
  \node{G_\omega}\arrow{ne,r,..}{p_B^{ab}}
\end{diagram} \]
The lemma also shows that if $g \in G_\omega$ can be represented by
some word in $\Ker(\rho_B^{ab})$ then, for any representation $g=F$
where $F$ is a word over $S$, $F$ is in $\Ker(\rho_B^{ab})$. We
identify $\Ker(\rho_B^{ab})$ with a subset of $G_\omega$ using
$\rho_\omega$, whence $\Ker(p_B^{ab})=\Ker(\rho_B^{ab})$ and $G_\omega
/\Ker(\rho_B^{ab}) \cong G_B^{ab}$.

\begin{proof}[Proof of lemma~\ref{lemma:welldef'd}]
  The proof is by induction on the length of $F$ and it will be done
  for all $\omega$ simultaneously.
  
  The statement is clear for the empty word. Next, no word of length
  $1$ represents the identity in any group $G_\omega$. Now assume that
  the claim is true for all words of length less than $n$, with
  $n\geq2$ and let $F$ be a word of length $n$ representing the
  identity in $G_\omega$.

  If $F$ is not reduced, then we reduce it to a shorter word $F'$.
  Since the reduced word $F'$ is in $\Ker(\rho_B^{ab})$ if and only if
  the original word $F$ is in $\Ker(\rho_B^{ab})$ the claim follows in
  this case from our inductive hypothesis.
 
  Assume $F$ is reduced. Since $F$ must be in $H_\omega$ the
  decomposition of $F$ of depth $1$ is well defined. The length of
  each of the (possibly not reduced) words
  $\overline{F_1},\dots,\overline{F_q}$ is at most $(n+1)/2<n$ (since
  $n>1$), so that, by the inductive hypothesis, each of the words
  $\overline{F_i}$ is in $\Ker(\rho_B^{ab})$.  The set
  $\Ker(\rho_B^{ab})$ is clearly closed under concatenation, so
  $\overline{F_1}\cdots\overline{F_q}$ is in $\Ker(\rho_B^{ab})$.
  Assume $F_B = b_1\cdots b_k$; then each of the $B$-letters $b_i$ from
  $F$ appears exactly once in some word $\overline{F_j}$, and
  therefore $F_B$ and $(\overline{F_1}\cdots \overline{F_q})_B$ 
  represent the same element in $G_B^{ab}$, namely, the identity.
  We conclude that $F$ is in $\Ker(\rho_B^{ab})$.
\end{proof} 

We may define $\Ker(\rho_A^{ab})$ and $p_A^{ab}:G_\omega \rightarrow
G_A^{ab}$ similarly. It is easy to see that any word representing the
identity in $G_\omega$ must come from $\Ker(\rho_A^{ab})$, so
$\Ker(p_A^{ab})=\Ker(\rho_A^{ab})$ and $G_\omega /\Ker(\rho_A^{ab}) \cong
G_A^{ab}$.

Since $G_\omega /\Ker(\rho_A^{ab}) \cong G_A^{ab}$ and $G_\omega
/\Ker(\rho_B^{ab}) \cong G_B^{ab}$ are abelian, the commutator subgroup
$[G_\omega,G_\omega]$ is in the intersection $\Ker(\rho_A^{ab}) \cap
\Ker(\rho_B^{ab})$. On the other hand, any word $F$ from the
intersection $\Ker(\rho_A^{ab}) \cap \Ker(\rho_B^{ab})$ clearly
represents the identity in the abelianization $G_\omega^{ab}$. We have
thus proved:

\begin{theorem} \label{theorem:commutators}
  For every $\omega$ in $\Omegah$, the commutator
  $[G_\omega,G_\omega]$ is equal to the intersection $\Ker(\rho_A^{ab})
  \cap \Ker(\rho_B^{ab})$. Moreover, $G_\omega^{ab} \cong G_A^{ab}
  \times G_B^{ab}$.
  
  As a consequence, the commutator $[G_\omega,G_\omega]$ is generated,
  as a subgroup of $G_\omega$, by all $[x,y]$ with $x,y\in G_A\cup
  G_B$.
\end{theorem} 

Let us define another set of words over $S$:
\[ \Ker(\rho_B) = \{ F\in S^*|\,F_B\text{ represents the identity in }G_B\}.\] 
Again, we can consider this set as a set of elements in $G_\omega$. It
is easy to see that this set is actually the normal closure
$G_A^{G_\omega}$ of $G_A$ in $G_\omega$ and it is generated as a
monoid by the set $\{a^g | a\in A, g\in G_B\}$. In case $G_B$ is
abelian, the sets $\Ker(\rho_B^{ab})$ and $\Ker(\rho_B)$ clearly
coincide and we have
\[G_\omega/G_A^{G_\omega} = G_\omega/\Ker(\rho_B) = G_\omega/\Ker(\rho_B^{ab}) 
    \cong G_B^{ab}=G_B.\] 
In particular, this shows that the index of the normal closure of $G_A$ is 
$p^2$ for any Grigorchuk $p$-group (as defined 
in~\cite{grigorchuk:gdegree} or~\cite{grigorchuk:pgps}). 

Of course $\Ker(\rho_A) = \{F\in S^*|\,F_A\text{ represents the
  identity in }G_A\}$ is the normal closure of $G_B$ in $G_\omega$,
but this is the subgroup $H_\omega$ which we already discussed.

The following result generalizes the decomposition of the first Grigorchuk 
group as a semidirect product. This approach allows a much more algebraic 
treatment of groups acting on trees. 
\begin{proposition}\label{TD}
  Let $G_A$ and $G_B$ be abelian groups and suppose $G_B$ splits as
  $G_B= K_1^\perp\times K_1$. Further, set $D=\langle A, K_1 \rangle$
  and $T=(K_1^\perp)^{G_\omega}$, the normal closure of
  $K_1^\perp$ in $G_\omega$. Then $G_\omega= T \rtimes D$.
  
  Moreover, the index of $T$ in $G_\omega$ is
  $|D|=|G_B|^q/|G_A|^{q-1}=|K_1|^q|G_A|$.
\end{proposition}
\begin{proof} First, note that $T^{G_\omega}$ is generated as a
  subgroup by the set
  \[ X=\{t^g,[t^g,d]|\, t\in K_1^\perp, g \in G_A, d \in D\}. \]
  Indeed, conjugation of any generator in $X$ by an element from $G_A$
  gives another generator, conjugation by $t\in K_1^\perp$ is
  unimportant since $K_1^\perp \subseteq X$ and, for $k \in K_1$,
  \begin{gather*} kt^gk^{-1} = [k,t^g]t^g,\\
    k[t^g,d]k^{-1} = kt^gd(t^{-1})^gd^{-1}k^{-1}=
    [k,t^g]t^gkd(t^{-1})^g(kd)^{-1} = [k,t^g][t^g,kd].
  \end{gather*} 
  The subgroup generated by $X$ is thus normal. On the other hand $t^g
  \in T$ and $[t^g,d]=t^g(t^{-1})^{dg} \in T$ so that $X \subseteq T$
  and $\langle X \rangle = T$.
  
  Since $G_\omega$ is generated by $D$ together with the elements of
  the form $t^g$, $t\in K_1^\perp$, $g\in G_A$, we note that
  $G_\omega= TD$.
  
  Let us prove that $T \cap D = 1$. Assume $g\in T\cap D$. Since $T
  \subseteq H_\omega$ we can consider $\psi(g)=(g_1,\dots,g_q)$. Since
  $g\in D$ we have $g_i \in K_1$, for all $i$, and
  $g_i=p_B^{ab}(g_i)$. On the other hand, since $g\in T$ we have
  $p_B^{ab}(g_i)\in \langle K_1^\perp \rangle$, for all $i$.
  Therefore $g_i=1$ for all $i$, so $g=1$.
  
  Consider $H_\omega \cap D$. Clearly $\psi(H_\omega \cap D) \subseteq
  K_1^q$. On the other hand, given $(k_1,\dots,k_q)\in K_1^q$ we have
  $\psi(k_1^{g_1}k_2^{g_2}\dots k_q^{g_q})=(k_1,\dots,k_q)$ where
  $g_i\in G_A$ with $g_i(q)=i$. Therefore $|H_\omega \cap D|=|K_1|^q$
  and since the index of $H_\omega$ in $G_\omega$ is $|G_A|$ we obtain
  the result.
\end{proof}

%****************************************************************
\part{Quantitative Estimates}
In the following sections we will impose various restrictions on the 
sequence $\omega$ defining the group $G_\omega$ and give estimates of word 
and period growth in those cases. All the estimates will be done with 
respect to the canonical generating set $S_\omega=A\cup B_\omega$. As a 
shorthand, we will use $\gamma_\omega(n)$ and $\pi_\omega(n)$ instead of 
$\gamma_{G_\omega}(n)$ and $\pi_{G_\omega}(n)$. 

%-----------------------------------------------------------------
\section{The Word Growth in the General Case}
A finite subsequence $\omega_{i+1}\omega_{i+2}\dots \omega_{i+r}$ of a
sequence $\omega$ in $\Omegah$ is \mathem{complete} if each element of
$G_B$ is sent to the identity by at least one homomorphism from the
sequence $\omega_{i+1}\omega_{i+2}\dots\omega_{i+r}$, i.e.\ if
$\bigcup_{j=1}^r \Ker(\omega_{i+j}) = G_B$. We note that a complete
sequence must have length at least $q+1$ since all the kernels have index 
$|G_A| \geq q$ in $G_B$. In particular, the length of a complete sequence 
is never shorter than $3$.  Note that by definition all sequences in 
$\Omegah$ can be factored into finite complete subsequences. 

\begin{theorem} \label{hatgrowth}
  $G_\omega$ has subexponential growth, for all $\omega$ in
  $\Omegah$.
\end{theorem}

Let $F$ be a reduced word of length $n$ representing an element in 
$H_\omega^{(r)}$ and consider the decomposition of the word $F$ of depth 
$r$. For $\ell=0,\dots,r$, define the length $|L_\ell(F)|$ of the level 
$\ell$ to be the sum of the lengths of the elements on the level $\ell$. 
The following lemma is a direct generalization of Lemma~1 
in~\cite{grigorchuk:pgps}. The proof is similar, but adapted to the more 
general setting of the present paper. 

\begin{lemma} [3/4-Shortening] \label{3/4} 
  Let $\omega \in \Omegah$ be a sequence that starts with a complete
  sequence of length $r$. Then the following inequality holds for
  every reduced word $F$ representing an element in $H^{(r)}_\omega$:
  \[ |L_r(F)| \leq  \frac{3}{4}|F| + q^r. \]
\end{lemma}
\begin{proof} 
  Define $\xi_i$ to be the number of $B$-letters from
  $K_i-(K_{i-1}\cup\dots\cup K_1)$ appearing in the words at the level
  $i-1$ and $\nu_i$ to be the number of simple reductions performed to
  get the words $F_{j_1\dots j_i}$ on the level $i$ from their
  unreduced versions $\overline{F_{j_1\dots j_i}}$.
  
  A reduced word $F$ of length $n$ has at most $(n+1)/2$ $B$-letters.
  Every $B$-letter in $F$ that is in $K_1$ contributes one
  $B$-letter and no $A$-letters to the unreduced words
  $\overline{F_1},\dots,\overline{F_q}$. The $B$-letters in $F$ that
  are not in $K_1$ (there are at most $(n+1)/2-\xi_1$ such letters) 
  contribute one $B$-letter and one $A$-letter.
  Finally, the $\nu_1$ simple reductions reduce the number of letters
  on level 1 by at least $\nu_1$. Therefore,
  \[ |L_1(F)| \leq 2((n+1)/2-\xi_1) + \xi_1 - \nu_1 = n+1-\xi_1-\nu_1. \]
  In the same manner, each of the $\xi_2$ $B$-letters on level $1$ that
  is from $K_2-K_1$ contributes at most one $B$-letter to the words on
  level $2$ and the other $B$-letters (at most $(|L_1(F)|+q)/2-\xi_2$ of
  them) contribute at most $2$ letters, so
  \[ |L_2(F)| \leq n+1+q-\xi_1-\xi_2-\nu_1-\nu_2. \]
  Proceeding in the same manner, we obtain the estimate 
  \begin{equation} \label{L}
    |L_r(F)| \leq n+1+q+\dots q^{r-1}- \xi_1-\xi_2-\dots-\xi_r-\nu_1-\nu_2-\dots-\nu_r. 
  \end{equation}
  If $\nu_1+\nu_2+\dots+\nu_r\geq n/4$, then the claim of the lemma follows. 
  Assume therefore
  \begin{equation} \label{nusmall}
    \nu_1+\nu_2+\dots+\nu_r < n/4.
  \end{equation}
  For $i=0,\dots,r-1$, define $|L_i(F)|^+$ to be the number of
  $B$-letters from $B-(K_1\cup\dots\cup K_i)$ appearing in the words
  at the level $i$. Clearly, $|L_0(F)|^+$ is the number of $B$-letters
  in $F$ and
  \[ |L_0(F)|^+ \geq \frac{n-1}{2}. \]
  Going from the level $0$ to the level $1$, each $B$-letter contributes
  one $B$ letter of the same type. Thus, the words
  $\overline{F_1},\dots,\overline{F_q}$ from the first level before the
  reduction takes place have exactly $|L_0(F)|^+ -\xi_1$ letters that
  come from $B-K_1$. Since we lose at most $2\nu_1$ letters due to the
  simple reductions, we obtain
  \[ |L_1(F)|^+ \geq \frac{n-1}{2} -\xi_1 - 2\nu_1. \]
  Next we go from level $1$ to level $2$. There are $|L_1(F)|^+$
  $B$-letters on level $1$ that come from $B-K_1$, so there are
  exactly $|L_1(F)|^+ -\xi_2$ $B$-letters from $B- (K_1 \cup K_2)$ in
  the words $\overline{F_{11}},\dots,\overline{F_{qq}}$ and then we
  lose at most $2\nu_2$ $B$-letters due to the reductions. We get
  \[ |L_2(F)|^+ \geq \frac{n-1}{2} -\xi_1 - \xi_2 - 2\nu_1 - 2\nu_2, \]
  and, by proceeding in a similar manner, 
  \begin{equation} \label{L+}
    |L_{r-1}(F)|^+ \geq \frac{n-1}{2} -\xi_1-\dots\xi_{r-1} -2\nu_1-\dots 2\nu_{r-1}. 
  \end{equation}
  Since $\omega_1\dots\omega_r$ is complete, we have $\xi_r=|L_{r-1}(F)|^+$ 
  and the inequalities~(\ref{L}),~(\ref{nusmall}) and~(\ref{L+}) give
  \[ |L_r(F)| \leq \frac{n}{2}+\frac{1}{2}+1+q+\dots q^{r-1} + 
  \nu_1 + \dots+\nu_{r-1} -\nu_r. \]
  which implies our claim. 
\end{proof}

We finish the proof of Theorem~\ref{hatgrowth} using either the
argument given in~\cite{grigorchuk:pgps} or the one
in~\cite[Theorem~VIII.61]{harpe:cgt}; namely let
\[e_\omega=\limsup_{n\to\infty}\sqrt[n]{\gamma_\omega(n)}\]
denote the \mathem{exponential growth rate} of $G_\omega$. It is known 
that this rate is 1 if and only if the group in question has 
subexponential growth. By the previous lemma we have $e_\omega\le 
e_{\sigma^r\omega}^{3/4}$ and since the $e_\omega$ are bounded (for 
instance, by $|A\cup B|$), it follows that $e_\omega=1$ for all 
$\omega\in\Omegah$. 

\subsection{A lower bound for word growth}
A general lower bound, tending to $e^n$ when $q\to\infty$, exists on
the word growth, and holds for all spinal groups:
\begin{theorem}\label{thm:wdlowerbd}
  $G_\omega$ has superpolynomial growth, for all $\omega$ in
  $\Omegah$. Moreover, the growth of $G_\omega$ satisfies
  \[e^{n^\alpha}\precsim\gamma_\omega(n),\]
  where $\alpha = \frac{\log(q)}{\log(q)-\log\frac12}$.
\end{theorem}

\begin{proof}
  Let $\gamma_\omega=\gamma_\omega^{(S,|\cdot|,\rho)}$ denote the growth
  of $G_\omega$ with respect to word length. We will obtain
  \begin{equation}\label{bigomega}
    \gamma_\omega(2qn+K_\omega)\ge L_\omega\gamma_{\sigma\omega}(n)^q\text{ for all }n\in\R_{\ge0}
  \end{equation}
  for some positive constants $K_\omega,L_\omega\ge0$.  Then, $x$ applications 
  of~(\ref{bigomega}) yield, neglecting the (unimportant)
  constant $K_\omega$, $\gamma_\omega((2q)^x)\ge
  L_\omega^{1+q+\dots+q^{x-1}}\gamma_{\sigma^x\omega}(1)^{q^x}$, from
  which the theorem's claim follows.
  See~\cite[Corollary~9]{bartholdi:lowerbd} for a similar proof.
  
  We now prove~(\ref{bigomega}). Choose some $h\in A$ with $h(1)=q$,
  let $\nu$ be a (set) retraction $A\to G_B$ of $\omega_1$, and
  consider the (set) map
  \[\lambda:\begin{cases}
    A\ni a&\mapsto h\nu(a)h^{-1}=\nu(a)^h\\
    B\ni b&\mapsto b.
  \end{cases}\]
  defined on reduced words over $A\cup B$. Note that $\lambda$ does
  not in general induce a group homomorphism, though this is the case
  for the first Grigorchuk group, where it is traditionally called
  $\sigma$. We may, however, naturally consider $\lambda(F)\in H_\omega$.
  
  Given any reduced word $F$ representing $x\in G_{\sigma\omega}$, we
  obtain an element $y=\lambda(F)$ of $H_\omega$, that has the
  following properties:
  \begin{align*}
    \varphi_q^\omega(y)&=x;\\
    \text{if $h(q)\neq1$, then }\varphi_1^\omega(y)&\in G_A\text{ and
    }\varphi_{h(q)}^\omega(y)\in\nu(W),\\
  &\text{ where $W$ is the set of $A$-letters in }x;\\
    \text{if $h(q)=1$, then }\varphi_1^\omega(y)&\in\langle A,\nu(W)\rangle;\\
    \varphi_i^\omega(y)&=1\text{ for all }i\not\in\{1,q,h(q)\}.
  \end{align*}
  
  In case $h(q)=1$, we restrict our consideration to words $F$ such
  that at most $q$ of their $A$-letters are not in $\langle h\rangle$.
  It then follows in all cases that all coordinates except the $q$-th
  of $\lambda(F)$ are bounded. To prove the only non-trivial case,
  suppose $h(q)=1$ and $F=F_0a_1F_1\dots a_qF_q$, where the $F_i$ are
  words over $\langle h\rangle\cup B$. Then
  \[\varphi_1^\omega(\lambda(F))=\alpha_0\nu(a_1)\alpha_1\dots\nu(a_q)\alpha_q,\]
  where $\alpha_i\in\langle A,\nu(h) \rangle$. By
  Corollary~\ref{finite-subgp}, this last group is finite, so each
  $\alpha_i$ is bounded (say of length at most $M$); then
  $\varphi_1^\omega(\lambda(F))$ is bounded, of length at most
  $N=q+(q+1)M$.
  
  Given words $F_1,\dots,F_q$ each of length at most $n$, we wish to
  construct a word $F$ with $\psi^\omega(F)=F'(F_1,\dots,F_q)F''$
  where $F'$ and $F''$ belong to a finite set $\mathcal F$. Let us
  choose elements $a_i$, for $i\in Y$, such that $a_i(q)=i$.  We take
  \[F=\lambda(F_1)^{a_1}\dots\lambda(F_{q-1})^{a_{q-1}}\lambda(F_q).\]
  In every coordinate $i$, we get $F_i$ (as
  $\varphi_i^\omega(\lambda(F_i)^{a_i})$), and other words, each of
  which is bounded. We may therefore take for $\mathcal F$ the set of
  words of length at most $(q-1)N$.

  Note now that for any word $F_i$ we have $|\lambda(F_i)|\le
  2|F_i|+1$, so $|F|\le q(2n+1)+q$. Also, $F$ has at most $q$
  $A$-letters not in $\langle h \rangle$, namely
  $a_1,a_1^{-1}a_2,\dots,a_{q-1}^{-1}$. Finally, $F$ determines
  $F_1,\dots,F_q$ up to the choice of $F'$ and $F''$, so~(\ref{bigomega})
  holds with $K=2q$ and $L=|\mathcal F|^{-2}$ 
\end{proof}

Note that if $q$ is a prime power $p^\alpha$, the group $G_\omega$ is
residually-$p$ so has growth at least $e^{\sqrt n}$
by~\cite{grigorchuk:hp}. The previous result is an improvement in
all cases but $q=2$.

For the special case of the first Grigorchuk group slightly better results 
exist, due to Yuri\u\i\ Leonov~\cite{leonov:pont} who obtained 
$\gamma(n)\succsim e^{n^{0.5041}}$, and to the first 
author~\cite{bartholdi:lowerbd} who obtained $\gamma(n)\succsim 
e^{n^{0.5157}}$. There is no doubt that a similar improvement of
Theorem~\ref{thm:wdlowerbd} for general spinal groups is possible.

%-----------------------------------------------------------------------------
\section{The Word Growth in the Case of Homogeneous Sequences}\label{sec:homo}
A sequence $\omega$ in $\Omegah$ is \mathem{$r$-homogeneous} (for
$r\geq3$) if every finite subsequence of length $r$ is complete. The set 
of $r$-homogeneous sequences in $\Omegah$ will be denoted by $\Omegar$.  
Note that $\Omegar$ is closed under the shift $\sigma$, a fact that is 
crucial for the arguments that follow, but we will not mention it 
explicitly anymore. Implicitly, all sequences $\omega$ in this section 
will come from $\Omegar$ for some fixed $r$.  We will prove the following: 

\begin{theorem} [$\eta$-Estimate] \label{(r)growth}
  If $\omega$ is an $r$-homogeneous sequence, then the growth function of
  the group $G_\omega$ satisfies
  \[ \gamma_\omega(n) \precsim e^{n^{\alpha}} \]
  where $\alpha= \frac{\log(q)}{\log(q)-log(\eta_r)}$ and $\eta_r$ is
  the positive root of the polynomial $x^r + x^{r-1} + x^{r-2} -2$. 
\end{theorem}

%-------------------------------------
\subsection{A triangular weight function on $G_\omega$}
The following weight assignment generalizes the approach taken
in~\cite{bartholdi:upperbd} by the first author in order to estimate
the growth of the first Grigorchuk group.

Consider the linear system of equations in the variables
$\tau_0,\dots,\tau_r$:
\begin{equation} \label{system}
  \begin{cases}
    \eta_r(\tau_0+\tau_i) = \tau_0+\tau_{i-1} & \text{ for }i=r,\dots,2,\\
    \eta_r(\tau_0+\tau_1) = \tau_r.
  \end{cases}
\end{equation}
The solution is given, up to a constant multiple, by
\begin{equation} \label{solution}
  \begin{cases}
    \tau_i = \eta_r^r + \eta_r^{r-i} -1 & \text{ for }i=r,\dots,1,\\
    \tau_0 = 1-\eta_r^r.
  \end{cases}
\end{equation}
If we also require $\tau_1+\tau_2=\tau_r$ we get that $\eta_r$ must be
a root of the polynomial $x^r+x^{r-1}+x^{r-2}-2$. If we choose
$\eta_r$ to be the root of this polynomial that is between $0$ and $1$
we obtain that the solution~(\ref{solution}) of the
system~(\ref{system}) satisfies the additional properties
\begin{align}
  0 < \tau_1 < \dots < \tau_r < 1, \qquad 0 < \tau_0 < 1, \\
  \tau_i+\tau_j \geq \tau_k \text{ for all } 1 \leq i,j,k \leq r\text
  { with }i\neq j.
  \label{triangular}
\end{align}
The index $r$ in $\eta_r$ will be sometimes omitted without warning.

Now, given $\omega \in \Omegar$, we define the weight of the generating 
elements in $S_\omega$ as follows: $\tau(a)=\tau_0$, for $a$ in $A$ and 
$\tau(b_\omega)=\tau_i$, where $i$ is the smallest index with 
$\omega_{i}(b)=1$. 

Clearly, $\tau$ is a triangular weight function. The only point worth
mentioning is that if $b$ and $c$ are two $B$-letters of the same
weight and $bc=d \not = 1$ then $d$ has no greater weight than $b$ or
$c$ (this holds because $b_\omega,c_\omega \in \Ker(\omega_i)$ implies
$d_\omega \in \Ker(\omega_i)$).

For obvious reasons, the weight
$\partial_{G_\omega}^{(S_\omega,\tau,\rho)}(g)$ for $g \in G_\omega$
will be denoted by $\partial^{\tau}(g)$ and, more often, just by
$\partial(g)$.

%------------------------------------
\subsection{A tree representation of the elements in $G_\omega$}
Let $g$ be an element in $G_\omega$. There is a unique element $h$ in
$G_A$ such that $hg$ is in $H_\omega$. We extend the map $\psi$ to
$G_\omega$ by $\psi(g)=\psi(hg)$ and we write
$\psi(g)=(g_1,\dots,g_q)$ in this case. This notation does not
interfere (too much) with our previous agreement since $h=1$ for $g$
in $H_\omega$. Note that the extended $\psi$ is not a homomorphism
(nor is it injective: it is $|G_A|$-to-one).

\begin{lemma} [$\eta$-Shortening] \label{<eta}
Let $g \in G_\omega$. Then
\[ \sum_{i=1}^q \partial^{\tau}(g_i) \leq \eta_r\big(\partial^{\tau}(g)+ \tau_0\big). \]
\end{lemma}
\begin{proof} Let a minimal form of $g$ be
  \[ g = [a_0]b_1a_1 \dots b_{k-1}a_{k-1}b_k[a_k]. \]
  Then $hg$ can be written in the form $hg = h[a_0]b_1a_1 \dots
  a_{k-1}b_k[a_k]$ and rewritten in the form
  \begin{equation} \label{hg} hg = b_1^{g_1} \dots b_k^{g_k},
  \end{equation}
  where $g_i=h[a_0]a_1 \dots a_{i-1} \in G_A$. Clearly, $\partial(g)
  \geq (k-1)\tau_0 +\sum_{j=1}^k \tau(b_j)$, which yields
  \begin{equation} \label{eta-interm}
    \sum_{i=1}^k \eta(\tau_0 + \tau(b_j)) \leq
    \eta(\partial(g)+\tau_0).
  \end{equation}
  Now, observe that if the $B$-generator $b$ is of weight $\tau_i$
  with $i>1$ then $\psi(b^g)$ has as components one $B$-generator of
  weight $\tau_{i-1}$ and one $A$-generator (of weight $\tau_0$ of
  course) with the rest of the components trivial. Thus, such a $b^g$
  (from~(\ref{hg})) contributes at most
  $\tau_0+\tau_{i-1}=\eta(\tau_0+\tau(b))$ to the sum $\sum
  \partial(g_i)$. On the other hand, if $b$ is a $B$-generator of
  weight $\tau_1$ then $\psi(b^g)$ has as components one $B$-generator
  of weight at most $\tau_r$, and the rest of the components are
  trivial. Such a $b^g$ contributes at most
  $\tau_r=\eta(\tau_0+\tau(b))$ to the sum $\sum \partial(g_i)$.
  Therefore
  \begin{equation} \label{eta-interm2}
    \sum_{i=1}^q \partial^{\tau}(g_i) \leq \sum_{j=1}^k \eta(\tau_0+\tau(b_j))
  \end{equation}
  and the claim of the lemma follows by combining~(\ref{eta-interm})
  and~(\ref{eta-interm2}).
\end{proof}

A simple corollary of the lemma above is that for any $\zeta$ with
$\eta<\zeta<1$ there exists a positive constant $K_\zeta=\eta
\tau_0/(\zeta-\eta)$ such that $\sum_{i=1}^q \partial(g_i) \leq \zeta
\partial(g)$ for every $g$ in $G_\omega$ with $\partial(g) \geq
K_\zeta$, and, therefore,
\begin{equation} \label{<zeta}
   \partial(g_i) \leq \zeta \partial(g)
\end{equation}
for all $i=1,\dots,q$ and every $g$ in $G_\omega$ with $\partial(g)
\geq K_\zeta$. For the discussion that follows just pick a fixed value
for $\zeta$.

Starting with an element $g$ in $G_\omega$ we construct a rooted,
$q$-regular, labeled tree $\iota(g)$ each of whose leaves is decorated
by an element of weight $\leq K$ (for a chosen $K \geq K_\zeta$) and
each of whose interior vertices is decorated by an element of $G_A$.
Technically, $\iota$ depends on $K$ but we choose not to indicate this
in our notation.

The tree $\iota(g)$ is called the \mathem{portrait} of $g$ of size $K$
and is constructed inductively as follows: if the weight of $g$ is
$\leq K$ then the portrait of $g$ is the tree that has one vertex
decorated by $g$; if $\partial(g)>K$, then the weight of all $g_i$'s
is at most $\zeta \partial(g)$ (see~(\ref{<zeta})) and the portrait of
$g$ is the tree that has $h$ at its root and the trees
$\iota(g_1),\dots,\iota(g_q)$ attached at the branches below the root.

The map $\iota$ sending each $g$ in $G_\omega$ to its portrait is
injective (see the corresponding proof in~\cite{bartholdi:upperbd}).
The main points are that $\psi$ is injective on $H_\omega$ and for
every $g$ in $G_\omega$ the element $h\in G_A$ such that $hg$ is in
$H_\omega$ is unique.

\begin{lemma} \label{Ln} There exists a positive constant $K$ such that
  \[L(n) \precsim n^\alpha\]
  with $\alpha=\frac{\log{(q)}}{\log(q)-\log(\eta_r)}$, where $L(n)$
  is the maximal possible number of leaves in the portrait of size $K$
  of an element of weight at most $n$.
\end{lemma}
\begin{proof}
Let $\kappa=\eta \tau_0/(q-\eta)$. Choose $K$ so that $K \geq 
\max\big(q^{1/\alpha}+\kappa,K_\zeta\big)$ is big enough in order that
\[j + \left(\frac{q}{q-j}\right)^{\alpha-1}
\left(n-\kappa +\frac{\kappa j}{\eta}\right)^\alpha \leq
(n-\kappa)^\alpha \] be satisfied for all $n>K$ and all
$j=0,\dots,q-1$. Such a choice is possible because
$(q/(q-j))^{\alpha-1}<1$ for $j=1,\dots,q-1$ and the two expressions
are equal for $j=0$.

Define a function $L'(n)$ on $\R_{>0}$ by
\[ L'(n) = \begin{cases} 1 & \text{ if }n \leq K,\\
  (n-\kappa)^\alpha & \text{ if }n > K.
\end{cases}\]
We prove, by induction on $n$, that $L(n)\leq L'(n)$. If the weight of $g$ 
is $\leq K$, the portrait has $1$ leaf and $L'(n)=1$. Otherwise, the 
portrait of $g$ is made up of those of $g_1,\dots,g_q$. Let the weights of 
these $q$ elements be $n_1,\dots,n_q$. By Lemma~\ref{<zeta} we have $n_i 
\leq \zeta n$, so by induction the number of leaves in the portrait of 
$g_i$ is at most $L'(n_i)$, $i=1,\dots,q$ and the number of leaves in the 
portrait of $g$ is, therefore, at most $\sum_{i=1}^q L'(n_i)$. 

Suppose that $j$ of the numbers $n_1,\dots,n_q$ are no greater than $K$ 
and the other $q-j$ are greater than $K$, where $0 \leq j \leq q-1$.  
Without loss of generality we may assume $n_1,\dots,n_j\leq K < 
n_{j+1},\dots,n_q$. Using Jensen's inequality, Lemma~\ref{<eta}, the fact 
that $\eta^\alpha=q^{\alpha-1}$ and $0 < \alpha < 1$, and direct 
calculation, we see that 
\[ \begin{aligned} 
  \sum_{i=1}^q L'(n_i) &= j + \sum_{i=j+1}^q (n_i-\kappa)^\alpha \leq 
  j + (q-j)\left(\frac{1}{q-j}\sum_{i=j+1}^q (n_i-\kappa)\right)^\alpha\\
  &= j + \left(\frac{q}{q-j}\right)^{\alpha-1}
  \left(n -\kappa + \frac{\kappa j}{\eta}\right)^\alpha \leq (n-\kappa)^\alpha = L'(n),
\end{aligned} \]
where the last used inequality holds by the choice of $K$. 

In case none of $n_1,\dots,n_q$ is greater than $K$, we have $\sum_{i=1}^q 
L'(n_i)=q$ which is no greater than $(n-\kappa)^\alpha = L'(n)$, again, by the 
choice of $K$. 
\end{proof}

Lemmata~\ref{easylemma} and~\ref{<eta} can be summed up in a
\begin{scholium}
  The portrait of an element $g$ of weight $n$ is a tree of sublinear
  size $n^\alpha$ and logarithmic depth $\log_2(n)$.
\end{scholium}

%----------------
\subsection{Proof of Theorem~\ref{(r)growth} ($\eta$-Estimate)}
The number of elements in $G_\omega$ of weight at most $n$, i.e.\ the
number of elements in $B_\omega(n)$ is equal to the number of trees in
$\iota(B_\omega(n))$. This number is bounded above by the number of
labeled, rooted, $q$-regular trees with at most $L(n)$ leaves where each 
of the leaves is decorated by an element of weight at most 
$K$ and each interior vertex by an element of $G_A$. 

The number $N(m)$ of labeled, rooted, $q$-regular trees with exactly $m$ 
interior vertices (and, therefore, exactly $(q-1)m+1$ leaves) is 
$\frac{1}{qm+1}\binom{qm+1}{m} \sim e^m$ 
(see~\cite[page~1033]{gessel-s:handbook}).  Thus the number $D(m)$ of such 
trees with at most $m$ interior vertices is also $\sim e^m$. A tree with 
at most $L(n)$ leaves has at most $I(n)=(L(n)-1)/(q-1) \sim n^\alpha$ 
interior vertices, so that the number of labeled trees we are interested 
in is $\sim e^{n^\alpha}$. 

The decoration of the interior vertices can be done in at most
$|G_A|^{I(n)} \sim e^{n^\alpha}$ ways.

Finally, note that different leaves of a tree representation live on
different levels and therefore in different groups
$G_{\sigma^t\omega}$.  But the number of elements of weight at most
$K$ is bounded above by a finite number, denote it by $\gamma_f(K)$,
for all groups on at most $|A \cup B|$ generators. Thus the decoration
of the leaves can be chosen in at most $\gamma_f(K)^{L(n)} \sim
e^{n^\alpha}$ ways.

Therefore $\gamma_\omega(n) \precsim e^{n^\alpha}$ and the proof is complete. 

%-------------------------------
\subsection{Word growth in the case of factorable sequences}
We say that a sequence in $\Omegah$ is \mathem{$r$-factorable} if it can 
be factored in complete subsequences of length at most $r$. The set of 
$r$-factorable sequences will be denoted by $\Omegabr$. Clearly, $\Omegar 
\subset \Omegabr\subset \Omega^{(2r-1)}$ and the inclusions are proper. An 
upper bound on the degree of word growth can thus be obtained from 
Theorem~\ref{(r)growth}, but we can do slightly better if we combine 
Lemma~\ref{3/4} with the idea of portrait of an element. 

\begin{theorem} [3/4-Estimate] \label{[r]growth3/4} 
  If $\omega$ is an $r$-factorable sequence, then the growth function of
  the group $G_\omega$ satisfies
  \[ \gamma_\omega(n) \precsim e^{n^{\alpha}} \]
  where $\alpha=\frac{\log(q^r)}{\log(q^r)-log(3/4)}=\frac{\log(q)}{\log(q)-
    \log\left(\sqrt[\uproot{2} r]{3/4}\right)}$.
\end{theorem} 
\begin{proof}
  Let $\omega$ be an $r$-factorable sequence, factored in complete
  words of lengths $r_1,r_2,r_3,\dots$, with all $r_i\leq r$. We can
  define a modification of the portrait of an element by requiring
  that whenever we ``blow up'' a leaf on the level $r_1+r_2+\dots +
  r_{i-1}$ because its size is too big, we expand it $r_i$ levels down
  (i.e.\ the original word is expanded $r_1$ levels down, the words on
  the level $r_1$ are expanded $r_2$ levels, the words on the level
  $r_1+r_2$ are expanded $r_3$ levels, etc.) and obtain $q^{r_i} \leq
  q^r$ new leaves. An analog of Lemma~\ref{3/4} still holds, i.e.\ if
  the length of $g$ is at most $n$ then sum of the lengths of the
  elements at the newly obtained $q^{r_i}$ leaves is at most
  $\frac34n+O(1)$ (indeed, the fact that we multiply by an element
  from $G_A$ here and there does not increase the $B$-length of the
  words in question, so that the sum of the lengths at the newly
  obtained leaves is still at most $\frac34n$ plus a constant).
  Proceeding as before completes the proof of the theorem.
\end{proof}

The 3/4-Estimate is obtained only for the class of Grigorchuk $p$-groups 
defined by $r$-homogeneous (not $r$-factorable as above) sequences by 
Roman Muchnik and Igor Pak in \cite{muchnik-p:growth} by different means.  
On the other hand, their approach gives slightly better results for the 
Grigorchuk 2-groups defined by $r$-homogeneous sequences except in the 
case $r=3$ in which case the estimates coincide. 

We can provide a small improvement in a special case that includes all 
Grigorchuk 2-groups. Namely, we are going to assume that $\omega$ is an 
$r$-factorable sequence such that each factor contains three homomorphisms 
whose kernels cover $G_B$. Note that this is possible only when $q=2$. 
Also, note that in case $q=2$ we must have $G_A=\Z/2\Z$ since that is the 
only group that acts transitively and faithfully on the two-element set 
$Y=\{1,2\}$. Since $G_B$ is a subdirect product of several copies of $G_A$ 
we must have $G_B=(\Z/2\Z)^d$ for some $d\geq 2$, i.e.\ $G_B$ is an 
elementary abelian 2-group. 

\begin{lemma} [2/3-Shortening] \label{2/3} 
  Let $q=2$ and $\omega \in \Omegah$ be a sequence such that there
  exist $3$ letters $\omega_k$, $\omega_\ell$ and $\omega_m$, $1\leq
  k<\ell<m \leq r$, with the property that $K_k \cup K_\ell \cup
  K_m=G_B$. Then the following inequality holds for every reduced word
  $F$ representing an element in $H^{(r)}_\omega$:
  \[ |L_r(F)| < \frac{2}{3}|F| + \frac{2}{3}+ 3q^r. \]
\end{lemma} 
\begin{proof}
  For $i=0,\dots,r$, denote by $|L_i|_A$, $|L_i|_B$ and $|L_i|_{K_j}$
  the number of $A$-letters, $B$-letters and $B$-letters from $K_j$,
  respectively, in the words on the level $i$ of the decomposition of
  depth $r$ of $F$.

  Clearly, $|L_0|_B \geq |L_1|_B \geq \dots \geq |L_r|_B$. Also,
  \[ |L_{i+1}|_A \leq |L_i|_B - |L_i|_{K_{i+1}}, \] 
  since every $B$-letter from the level $i$ contributes at most one
  $A$-letter to the next level, except for those $B$-letters that are
  in $K_{i+1}$. This gives
  \[ |L_{i+1}|_B \leq |L_{i+1}|_A + q^{i+1} \leq 
  |L_i|_B - |L_i|_{K_{i+1}} + q^{i+1}, \] 
  and therefore
  \begin{align}
    |L_r|_B &\leq |L_{k-1}|_B - |L_{k-1}|_{K_k} + q^k - |L_{\ell-1}|_{K_\ell} +
    q^\ell - |L_{m-1}|_{K_m} + q^m \notag\\
    &< (|L_{k-1}|_{K_\ell} - |L_{\ell-1}|_{K_\ell}) + (|L_{k-1}|_{K_m}) - 
    |L_{m-1}|_{K_m} + 3q^r, \label{LrB<}
  \end{align}
  since $|L_{k-1}|_B \leq |L_{k-1}|_{K_k}+|L_{k-1}|_{K_\ell} + |L_{k-1}|_{K_m}$.  

  It is easy to see that 
  \[ |L_i|_{K_j} - |L_{i+1}|_{K_j} \leq |L_i|_B - |L_{i+1}|_B. \]
  Indeed, any change (up or down) in the number of letters in $K_j$
  going from the level $i$ to the level $i+1$ is due to simple
  reductions involving letters from $K_j$, but each such simple
  reduction also changes (always down) the total number of $B$-letters
  by the same amount. Then, by telescoping,
  \begin{equation} \label{<B}
    |L_i|_{K_j} - |L_t|_{K_j} \leq |L_i|_B - |L_t|_B, 
  \end{equation}
  whenever $i \leq t$. 

  Combining~(\ref{LrB<}) and~(\ref{<B}) gives
  \begin{align*}
    3|L_r|_B &\leq |L_r|_B+|L_{\ell-1}|_B+|L_{m-1}|_B\\
    &< (|L_{k-1}|_B - |L_{\ell-1}|_B) + (|L_{k-1}|_B -
    |L_{m-1}|_B) + 3q^r+|L_{\ell-1}|_B+|L_{m-1}|_B\\
    &\leq 2|L_{k-1}|_B+ 3q^r \leq 2|L_0|_B+3q^r,
  \end{align*}
  which implies our result since $|L_r(F)|\leq
  2|L_r|_B+q^r$ and $|L_0|_B \leq (|F|+1)/2$.
\end{proof}

Let us mention here that, for every $k$, there exists a reduced word $F$ 
of length $n=24k$ in the first Grigorchuk group representing an element in 
$H^{(3)}_\omega$ such that $|L_3(F)| = 16k = 2n/3$. An example of such a 
word is $(abadac)^{4k}$. Thus, the lemma above cannot be improved in the 
sense that there cannot be an improvement unless one starts paying 
attention to reductions that are not simple. Of course, $(abadac)^{16}=1$ 
in the first Grigorchuk group, so by introducing other relations the
multiplicative constant of Lemma~\ref{2/3} could be sharpened. 

As a corollary to the shortening lemma above, we obtain: 

\begin{theorem} [2/3-Estimate] \label{[r]growth2/3}
  Let $q=2$ and $\omega$ be an $r$-factorable sequence such that 
  each factor contains three letters whose kernels cover $G_B$. Then the
  growth function of the group $G_\omega$ satisfies
  \[ \gamma_\omega(n) \precsim e^{n^{\alpha}} \]
  where $\alpha=
  \frac{\log(q^r)}{\log(q^r)-log(2/3)}=\frac{\log(q)}{\log(q)-
    \log\left(\sqrt[\uproot{2} r]{2/3}\right)}$.
\end{theorem} 

The previous theorem and the lemma just before could be generalized to the 
other values of $q$, but the shortest possible complete sequence in those 
cases would be at least 4 and at best we would obtain the 3/4-Estimate 
already provided before. 

%-------------------------------
\subsection{Calculations and comparisons}
The Table \ref{comparison} included below lists various values of $\alpha$ 
(always rounded up), for different $q$ and $r$, such that 
$\gamma_\omega(n)\precsim e^{n^\alpha}$. The entries in the last column 
indicate conditions on $\omega$. A row labeled by ``homo.'' indicates that 
the estimate is valid for $r$-homogeneous sequences.  The estimate in a 
row labeled by ``fact.'' is valid for $r$-factorable sequences. In 
addition, the labels ``$\eta$'', ``3/4'' and ``2/3'' indicate which 
estimate is used. 

\begin{table} 
\begin{tabular}{|C|CCCCCCCC|c|}
    \hline
    q   &r=3  &r=4  &r=5  &r=6  &r=7  &r=8  &r=9  &r=10 & Condition\\
    \hline
    2   &.768 &.836 &.872 &.896 &.912 &.924 &.933 &.940 &homo. ($\eta$)\\ 
    2   &.837 &.873 &.896 &.912 &.923 &.932 &.939 &.945 &fact. (2/3)\\
    2   &.879 &.906 &.924 &.936 &.945 &.951 &.956 &.961 &fact. (3/4)\\
    \hline
    3   &     &.890 &.916 &.932 &.943 &.951 &.957 &.961 &homo. ($\eta$)\\
    3   &     &.939 &.951 &.959 &.964 &.969 &.972 &.975 &fact. (3/4)\\
    \hline
    4   &     &     &.932 &.945 &.954 &.960 &.965 &.969 &homo. ($\eta$)\\
    4   &     &     &.960 &.967 &.972 &.975 &.977 &.980 &fact. (3/4)\\
    \hline
    5   &     &     &     &.952 &.960 &.966 &.970 &.973 &homo. ($\eta$)\\
    5   &     &     &     &.972 &.976 &.979 &.981 &.983 &fact. (3/4)\\
    \hline
\end{tabular} \vspace{3mm}
\caption{Comparison of the obtained estimates} \label{comparison} 
\end{table}

Note that the $\eta$-estimate for an $r$-homogeneous sequence is always 
better than the $2/3$-Estimate or $3/4$-Estimate for an $r$-factorable 
sequence (since $\eta_r < \sqrt[\uproot{2} r]{2/3} 
<\sqrt[\uproot{2} r]{3/4}$).

Also, note that the best available $3/4$-Estimate is sometimes better than 
the best available $\eta$-Estimate. There are sequences $\omega$ is in 
$\Omegabr$ with the property that the smallest value of $s$ such that 
$\omega$ is $s$-homogeneous is $2r-1$. The $3/4$-Estimate for 
$r$-factorable sequences is better than the $\eta$-Estimate for 
$(2r-1)$-homogeneous sequences. For the sake of an example, take $q=3$ and 
$r=4$ and let $\omega=0123321001233210\dots$ where all homomorphisms 
$0,1,2,3$ are required to form a complete sequence. Then, the smallest 
value of $s$ such that $\omega$ is $s$-homogeneous is $s=7$ and the 
$3/4$-Estimate gives $\alpha=.939$ while the $\eta$-Estimate gives 
$\alpha=.943$. 

%-----------------------------------------------------------------------------
\section{Periodicity and Period Growth in case of a Regular
  Action}\label{sec:period} 
>From now on we assume, in addition to the earlier conditions, that the
action of $G_A$ on $Y$ is regular, i.e.\ each permutation of $Y$ induced 
by an element in $G_A$ is regular (recall a permutation is 
\mathem{regular} if all of its cycles have the same length). The number of 
elements in $G_A$ must be equal to $|Y|=q$ in case of a faithful and 
regular action. Also, note that if $G_A$ is abelian it must act regularly 
on $Y$.  

The group order of an element $g$ will be denoted by $\pi(g)$. In case
$F$ is a word $\pi(F)$ denotes the order of the element represented by
the word $F$.

%----------------
\subsection{Period decomposition and periodicity}
We will describe a step in a procedure introduced
in~\cite{grigorchuk:gdegree} that will help us to determine some
upper bounds on the period growth of the constructed groups.

Let $F$ be a reduced word of even length of the form  
\begin{equation} \label{crform}
F= b_1a_1\dots b_ka_k. 
\end{equation}

Rewrite $F$ in the form $F= b_1b_2^{g_2}\dots b_k^{g_k}a_1\dots a_k$,
where $g_i=a_1\dots a_{i-1}$, $i=2,\dots,k$. Set $g=a_1\dots a_k\in
G_A$ and let its order be $s$, a divisor of $q$. Note that $g=1$
corresponds to $F \in H_\omega$. The length of each cycle of $g$ is
$s$ because of the regularity of the action. Put $H=b_1b_2^{g_2}\dots
b_k^{g_k}$ and consider the element $F^s=(Hg)^s\in H_\omega$. We
rewrite this element in the form $F^s = HH^g\dots H^{g^{s-1}}$ and
then in the form
\begin{equation} \label{F^s}
F^s= \left(b_1b_2^{g_2}\dots b_k^{g_k}\right) \left(b_1^g(b_2^{g_2})^g\dots 
(b_k^{g_k})^g\right) \dots \left(b_1^{g^{s-1}}(b_2^{g_2})^{g^{s-1}}\dots 
(b_k^{g_k})^{g^{s-1}}\right). 
\end{equation} 
Next, by using tables similar to Table~\ref{table:phi} (but for all 
possible $a$), we calculate the (possibly unreduced) words 
$\overline{F_1},\dots,\overline{F_q}$ representing 
$\varphi_1(F^s)$,$\dots$, $\varphi_q(F^s)$, respectively, and then we use 
simple reductions to get the reduced words $F_1,\dots,F_q$. The tree that 
has $F$ at its root and $F_1,\dots,F_q$ as leaves on the first level is 
called the \mathem{period 
  decomposition} of $F$. Clearly $\psi(F^s)= (F_1,\dots,F_q)$ holds
and the order $\pi(F)$ of $F$ is divisor of $q \cdot
\gcd(\pi(F_1),\dots,\pi(F_q))$ since $\psi$ is injective, $s$ divides $q$ 
and the elements $(F_1,1,\dots,1), \dots, (1,\dots,1,F_q)$ commute in 
$\Pi_{i=1}^q G_{\sigma\omega}$. 

The notation introduced above for the vertices in the period decomposition 
interferes with the notation introduced before for the vertices in the 
decomposition of words, but we are not going to use the latter anymore. 

Let us make a couple of simple observations on the structure of the
possibly unreduced words $\overline{F_1},\dots,\overline{F_q}$ used to
obtain the reduced words $F_1,\dots,F_q$ of the period decomposition.

The conjugate elements $b_1,b_1^g,\dots,b_1^{g^{s-1}}$ appear in the
expression~(\ref{F^s}). The generator $b_1$ contributes exactly one
appearance of the letter $b_1$ to $\overline{F_q}$. The other
conjugates $b_1^g,\dots,b_1^{g^{s-1}}$ of $b_1$ contribute exactly one
appearance of the letter $b_1$ to the words
$\overline{F_{g(q)}},\dots,\overline{F_{g^{s-1}(q)}}$, respectively.
Similarly, $b_i^{g_i}$ contributes exactly one appearance of the
letter $b_i$ to the word $\overline{F_{g_i(q)}}$ and each of its
conjugates $(b_i^{g_i})^g,\dots,(b_i^{g_i})^{g^{s-1}}$ contributes
exactly one appearance of $b_i$ to the words
$\overline{F_{gg_i(q)}},\dots,\overline{F_{g^{s-1}g_i(q)}}$. Since the
length of each $g$-orbit is $s$ we see that, as far as the $B$-letters
are concerned, no word $\overline{F_i}$ gets more than one of each of
the letters $b_1,\dots,b_k$, possibly not in that order. Similarly, no
word $\overline{F_i}$ can get more than $k$ $A$-letters and it is
possible to get $k$ $A$-letters only if none of the letters
$b_1,\dots,b_k$ is in $K_1$. More precisely, the maximal number of
$A$-letters in any $\overline{F_i}$ is $k-|F|_{K_1}$.

\begin{theorem}
  Let $\omega$ be a sequence in $\Omegah$. Then the group $G_\omega$
  is periodic.
\end{theorem}
\begin{proof}
  We will prove that the order of any element $g$ in $G_\omega$
  divides some power of $q$.  The proof is by induction on the
  length $n$ of $g$ and it will be done for all $\omega$
  simultaneously.

  The statement is clear for $n=0$ and $n=1$. Assume that it is true
  for all words of length less than $n$, where $n\geq 2$, and consider
  an element $g$ of length $n$.
 
  If $n$ is odd the element $g$ is conjugate to an element of smaller
  length and we are done by the inductive hypothesis. Assume then that
  $n$ is even.  Clearly, $g$ is conjugate to an element that can be
  represented by a word of the form
  \[ F= b_1a_1\dots b_ka_k. \] 
  In this case $\pi(g)=\pi(F)$ divides
  $q\cdot \gcd(\pi(F_1),\dots,\pi(F_q))$ and if all the words $F_i$ have
  length shorter than $n$ we are done by the inductive hypothesis.

  Assume that some of the words $F_i$ have length $n$. This is
  possible only when $F$ does not have any $B$-letters from $K_1$.
  Also, the words $\overline{F_i}$ corresponding to the words $F_i$ of
  length $n$ must be reduced, so that the words $F_i$ having length
  $n$ have the same $B$-letters as $F$ does. For each of these
  finitely many words we repeat the discussion above; namely, for each
  such $F_i$ of length $n$ we construct the period decomposition.
  Either all of the constructed words $F_{ij}$ are strictly shorter than $n$,
  and we get the result by induction; or some have length $n$, but the
  $B$-letters appearing in them do not come from $K_1\cup K_2$.

  This procedure cannot go on forever since $K_1\cup K_2 \cup \dots
  K_r=G_B$ holds for some $r\in\N$.  Therefore at some stage we get a
  shortening in all the words and we conclude that the order of $F$ is
  a divisor of some power of $q$.
\end{proof}

% -----------------------------
\subsection{Period shadow and period growth in case of homogeneous sequences}
We can give a polynomial upper bound on the period growth of
$G_\omega$ in case $\omega$ is a homogeneous sequence. In order to do
so, we will make another use of the triangular weight function $\tau$
introduced before.

\begin{lemma}
  Let $F=b_1a_1b_2\dots b_ka_k$ be a reduced word of length $2k$. Then
  \begin{equation} \label{<etaperiod}  
    \tau(F_i) \leq \eta_r\tau(F), \quad \text{for all} \quad 1\leq i \leq q. 
  \end{equation}
\end{lemma}
\begin{proof}
  $\tau(F)= \sum_{i=1}^k (\tau_0 + \tau(b_i))$, yielding $\sum_{i=1}^k
  \eta(\tau_0 + \tau(b_i)) = \eta \tau(F)$. Using an argument similar
  to that in Lemma~\ref{<eta} and the observations on the structure of
  the words $\overline{F_i}$ given above, we conclude that
  \[\tau(F_i) \leq \sum_{i=1}^k \eta(\tau_0 + \tau(b_i)) = \eta \tau(F).\] 
\end{proof}

Note that all the canonical generators of $G_\omega$ have weight no
more than $1$. Given an element $g$ in $G_\omega$ we construct a
rooted, $q$-regular, labeled tree, whose leaves are decorated by
elements of weight at most $1$ and whose interior vertices are
decorated by divisors of $q$. We call such a tree $\shadow(g)$ a
\mathem{period shadow} of $g$ (of size $1$). Note that $\shadow(g)$ is
not uniquely determined by $g$ --- nor does it uniquely determine $g$.

A period shadow of $g$ is constructed inductively as follows: let $g'$
be an element of minimal weight in the conjugacy class of $g$. If
$\tau(g')\le1$ then the shadow is the tree with one vertex decorated
by $g'$; if $\tau(g')>1$, we assume $g'$ is represented by a word $F$
in the form~(\ref{crform}), from which we construct words
$F_1,\dots,F_q$ each of weight at most $\eta\partial(F)$ and group
order $s$ (the order of $a_1\dots a_k$) dividing $q$. A shadow of $g$
is the tree with $s$ at its root and $\shadow(F_1),\dots,\shadow(F_q)$
attached to the root.

Let $C$ be the $\gcd$ of the periods of all the elements of $G_\omega$
with weight at most $1$. If $\shadow(g)$ is a shadow of $g$, then
\begin{equation} \label{pi divides}
\pi(g) \text{ divides } Cq^d,
\end{equation} 
where $d$ is the \mathem{depth} of $\shadow(g)$, i.e.\ the length of
the longest path from the root to a leaf. If $g$ is an element of
$G_\omega$ of weight $n$, the depth of the shadow of $g$ cannot be
greater than $\lceil\log_{1/\eta}(n)\rceil$, so the following theorem
holds:

\begin{theorem}[Period $\eta$-Estimate] \label{(r)pgrowth} 
  If $\omega$ is an $r$-homogeneous sequence, then the period growth
  function of the group $G_\omega$ satisfies
  \[ \pi_\omega(n) \precsim n^{\log_{1/\eta_r}(q)} \]
  where $\eta_r$ is the positive root of the polynomial $x^r + x^{r-1}
  + x^{r-2} -2$.
\end{theorem}

In a similar manner we can prove the following two theorems. 

\begin{theorem}[Period 3/4-Estimate]
  If $\omega$ is an $r$-factorable sequence, then the period \linebreak 
  growth function of the group $G_\omega$ satisfies
  \[ \pi_\omega(n) \precsim n^{r\log_{4/3}(q)}. \]
\end{theorem} 

Instead of a proof, let us just note that in the process of building a
shadow of size $1$ of an element $g$ of ordinary length $n$ we are
not sure that there is a shortening in the length at each level, but
there is a shortening by at least a factor of $3/4$ after no more than
$r$ levels. Thus, the depth of such a shadow cannot be greater than
$r\lceil\log_{4/3}(n)\rceil$ and the claim follows.

\begin{theorem}[Period 2/3-Estimate] 
  If $q=2$ and $\omega$ is an $r$-factorable sequence such that 
  each factor contains three letters whose kernels cover $G_B$, 
  then the period growth function of the group $G_\omega$ satisfies
  \[ \pi_\omega(n) \precsim n^{r\log_{3/2}(q)}. \]
\end{theorem}

%------------------------------------
\subsection{Period growth in the case of a prime degree} 
In addition to the regularity requirement we assume that the degree $q$ of 
the tree $\tree$ is a prime number. Thus the root group $G_A$ is cyclic of 
prime order $q=p$ and there is no loss in generality if we assume that 
$G_A$ is generated by the cyclic permutation $a=(12\dots p)$. We assume 
all this in this subsection without further notice. 

Let us describe the construction of a sequence that we call the
\mathem{period sequence} of an element $g$ in $G_\omega$.

First we represent $g$ by a reduced word $F_g$. Then we conjugate
$F_g$ until we get either a word $F$ of length 1 in which case we
stop, or we get a cyclically reduced word $F$ of the
form~(\ref{crform}). This word either represents an element in
$H_\omega$ in which case we stop or it has to be raised to the $p$-th
power to get an element in $H_\omega$. Consider the latter case and
take a look again at the expression~(\ref{F^s}). Clearly,
\[ \varphi_1(F^p) = \varphi_1(H) \varphi_1(H^g) \dots \varphi_1(H^{g^{p-1}})= 
    \varphi_1(H) \varphi_{g^{-1}(1)}(H) \dots \varphi_{g^{1-p(1)}}(H), \]
where $H=b_1b_2^{g_2}\dots b_k^{g_k} \in H_\omega$. Similarly, 
\[ \varphi_i(F^p) = \varphi_i(H) \varphi_i(H^g) \dots \varphi_i(H^{g^{p-1}})= 
    \varphi_i(H) \varphi_{g^{-1}(i)}(H) \dots \varphi_{g^{1-p(i)}}(H), \] 
so that all the elements $F_i$ from the period decomposition are conjugate and 
we have $\pi(g)=\pi(F)=p\pi(F_1)$. 

Each of the letters $b_1,\dots,b_k$ appears in $\overline{F_1}$. Also,
each of the $A$-letters or identity factors
$\omega_1(b_1),\dots,\omega_1(b_k)$ appears in $\overline{F_1}$. Thus, if 
$b=b_1b_2\dots b_k$ in $G_B$ we have $b_\omega=\rho_B(F)$ and 
$\rho_B(\overline{F_1})=\rho_B(F_1)=b_{\sigma\omega}$, because $G_B$ is 
commutative. We may also write 
$b=\rho_B(F)=\rho_B(\overline{F_1})=\rho_B(F_1)$, by dropping the indices 
as usual. 

On the other hand, we have
$\rho_A(\overline{F_1})=\rho_A(F_1)=\omega_1(b_1)\dots
\omega_1(b_k)=\omega_1(b)$.

We conjugate the word $F_1$ until we get a word of length $1$ or a
word of the form~(\ref{crform}). The conjugation does not change the
projections $\rho_B(F_1)=b$ and $\rho_A(F_1)=\omega_1(b)$ and can only
decrease the length.  Now, the cyclically reduced version of $F_1$ has
the same order as $F_1$ and either represents an element in $H_\omega$
or it has to be raised to the $p$-th power to get an element in
$H_\omega$. In the first case we stop. In the latter case we
construct, as before, a word $F_{11}$ such that
$\pi(g)=\pi(F)=p\pi(F_1)=p^2\pi(F_{11})$, $\rho_B(F_{11})=b$ and
$\rho_A(F_{11})=\omega_2(b)$.

This process cannot last forever, since $\omega_i(b)=1$ holds for some
$i$. The sequence $F,F_1,\dots,F_{\underbrace{1\dots 1}_{t}}$ obtained
this way has the property that the last word in the sequence, denoted
$F'$, has length $1$ or represents an element in $H_{\sigma^t\omega}$.
Also $\pi(F)=p^t\pi(F')$.

\begin{theorem} \label{r-1/1}
Let $q=p$ be a prime and $\omega$ an $r$-homogeneous word. 

If $p\geq 3$ or $p=2$ and each subsequence of $\omega$ of length $r$ 
contains three homomorphisms whose kernels cover $G_B$, then the period 
growth function of the regular spinal group $G_\omega$ satisfies 
  \[ \pi_\omega(n) \precsim n^{(r-1)\log_2(p)}. \]
\end{theorem} 
\begin{proof} 
  As usual, we use induction on $n$ and we prove the statement
  simultaneously for all $r$-homogeneous $\omega$. We will prove that
  $\pi_\omega(n) \leq Cp^{(r-1)\log_2(n)}$ where $C=p^2$.

  The statement is obvious for $n=1$. 
 
  Consider an element $g$ of length $n$, $n \geq 2$ and let
  $F,F_1,\dots,F_{\underbrace{1\dots 1}_{t}}=F'$ be its period sequence.
  We know that $t \leq r$ because the word $\omega$ is
  $r$-homogeneous.
 
  If $F'$ has length $1$ then it has order $p$ and $\pi(g)=p^{t+1} \leq
  p^{r+1}=Cp^{r-1} \leq Cp^{(r-1)\log_2(n)}$.
 
  Consider the case when $F'$ has (even) length greater than $1$ and
  $F_t \in H_{\sigma^t\omega}$. In that case $\pi(F') = \pi(W)$ for
  some $W$ in $G_{\sigma^{t+1}\omega}$ that has length at most half
  the length of $F'$ (the word $W$ is one of the leaves of the period
  decomposition of $F'$).

  In case $t \leq r-1$ we have $\pi(g)\leq
  p^{r-1}\pi_{\sigma^r\omega}(n/2)$ which is no greater than
  $Cp^{(r-1)\log_2(n)}$ by the induction hypothesis.
 
  Let $t=r$ and $p \geq 3$. Since $t=r$ we know that the length of
  $F'$ is at most $3n/4$ so that the length of $W$ is at most $3n/8$
  and we have $\pi(g) \leq p^r \pi_{\sigma^{r+1}\omega}(3n/8)$ which
  is no greater than $Cp^{(r-1)\log_2(n)}$ by the induction hypothesis
  and the fact that $r\geq 4$ holds in this case.
 
  In case $t=r$ and $p=2$, the length of $F'$ is at most $2n/3$ and we
  have $\pi(g) \leq p^r \pi_{\sigma^{r+1}\omega}(n/3)$ which is no
  greater than $Cp^{(r-1)\log_2(n)}$ by the induction hypothesis and
  the fact that $r \geq 3$ holds in this case.

  Therefore, in each case $\pi(g) \leq Cp^{(r-1)\log_2(n)} =
  Cn^{(r-1)\log_2(p)}$, which proves our claim.
\end{proof}

Note that the theorem above gives the estimate $\pi_\omega(n) \precsim 
n^{r-1}$ in case $q=p=2$ and every subsequence of $\omega$ of length $r$ 
contains 3 homomorphisms whose kernels cover $G_B$. In case this last 
condition does not hold we can still give the estimate $\pi_\omega(n) 
\precsim n^r$.  

%----------------------------------------------------------------
\subsection{Period growth for Grigorchuk $2$-groups}
We give here a tighter upper bound on the period growth of the Grigorchuk 
2-groups. It is based on a more precise observation of the process 
described in the previous subsection. 

\begin{theorem} \label{thm:period ub for 2-groups}
  Let $G_\omega$ be a Grigorchuk 2-group. If $\omega$ is an
  $r$-homogeneous word, then the period growth function of the group
  $G_\omega$ satisfies
  \[ \pi_\omega(n) \precsim n^{r/2}. \]
\end{theorem} 
\begin{proof}
Let $\chi_\omega:G_\omega \to G_\omega^{ab}$ be the abelianization map. 
Recall that $G_\omega^{ab}=\langle a\rangle\times\langle b,c\rangle$ is 
the elementary $2$-group of rank $3$. We recast the construction of the 
period sequence as follows: in the graph below, nodes correspond to images 
of elements $g$ under $\chi_*$; arrows indicate taking a projection,  
$\varphi_1$ or $\varphi_2$. Double arrows indicate a squaring was applied 
before taking the projection (because $g$ was not yet in $H_*$). Also, 
recall that in the squaring case the obtained projections $\varphi_1(g^2)$ 
and $\varphi_2(g^2)$ are conjugate. A condition labeling an edge indicates 
that such an edge can exist only if the condition is satisfied. 
\vspace{4mm} 

\newcommand\msmash[1]{\raisebox{-1ex}{\makebox(0,0)[b]{\smash{$#1$}}}}
\begin{center}
  \begin{picture}(225,80)(0,-10)
    \put(0,0){\msmash{(a,b)}} \put(13,-4){\vector(1,0){24}}\put(13,-1){\vector(1,0){24}}
    \put(0,20){\arc{25}{2.07}{7.35}\arc{31}{2.07}{7.35}}
    \put(-5,5){\vector(2,-1){0}}
    \put(-4,8){\vector(2,-1){0}}
    \put(12,-15){\footnotesize{$b \in K_t$}}
    \put(50,0){\msmash{(1,b)}} 
       \put(63,-2){\vector(1,0){24}} \put(50,8){\vector(1,2){16}} \put(50,8){\vector(-1,2){16}}
       \put(60,-15){\footnotesize{$b \not \in K_{t+1}$}} \put(60,19){\footnotesize{$b \in K_{t+1}$}} 
    \put(100,0){\msmash{(a,1)}} \put(113,-4){\vector(1,0){24}} \put(113,-1){\vector(1,0){24}}
        \put(112,-15){\footnotesize{$\bar{b} \in K_{t+2}$}}
    \put(150,0){\msmash{(1,1)}} 
       \put(163,-2){\vector(1,0){24}} \put(150,8){\vector(1,2){16}} \put(150,8){\vector(-1,2){16}}
       \put(160,-15){\footnotesize{$\bar{b} \not \in K_{t+3}$}} \put(160,19){\footnotesize{$\bar{b} \in K_{t+3}$}}
    \put(200,0){\msmash{(a,\bar{b})}} 
    \put(175,50){\msmash{(1,\bar{b})}} \put(175,58){\vector(1,2){5}}
    \put(125,50){\msmash{(1,1)}} \put(125,58){\vector(-1,2){5}}
    \put(75,50){\msmash{(1,1)}} \put(75,58){\vector(1,2){5}}
    \put(25,50){\msmash{(1,b)}} \put(25,58){\vector(-1,2){5}}
    \end{picture}
\end{center}
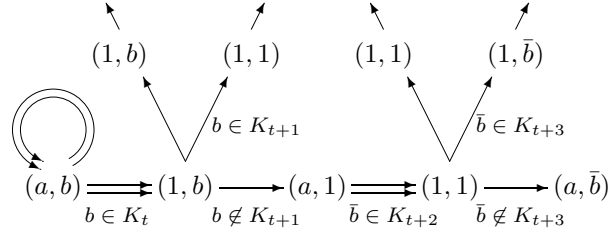
\centerline{One step of the algorithm for calculating the order of an 
element}\vspace{5mm} 
  
We proved in the previous section that all double arrows are as
described.  Let us complete the proof for single arrows.  According to
Theorem~\ref{theorem:commutators}, the commutator of a Grigorchuk
group is
\[[G_\omega,G_\omega] = \langle x=[a,b], y=[a,c], z=[a,d] \rangle.\]

Let us take an arbitrary element $g\in G_\omega$ with $\chi(g)=(a,b)$ and 
assume that $b$ is in the kernel $K_t$, but it is not in any kernel with 
smaller index. We follow our squaring procedure $t-1$ times obtaining an 
element $f \in G_{\sigma^{t-1}\omega}$ and set $h=\varphi_2(f^2)$. We know 
that $\pi(g)=2^t\pi(h)$, $\chi(h)=(1,b)$ and we wish to compute 
$\chi(\varphi_i(h))$ for $i=1,2$. For this purpose, write $f=u_1\dots u_m 
ab$ for some $u_i\in\{x,y,z\}$. Then   
\[h=\varphi(f^2)=\varphi_2(u_1)\dots\varphi_2(u_m)\varphi_2(u_1^a)\dots\varphi_2(u_m^a)b.\] 
Now, note that $x^a=x^{-1}$, $y^a=y^{-1}$, $z^a=z^{-1}$ and 
$\varphi_2(x)=\varphi_2(x^a)=b$, $\varphi_2(y)=ca$, $\varphi_2(x^a)=ac$, 
$\varphi_2(y)=da$, $\varphi_2(x^a)=ad$, so that all $B$-letters appear in 
pairs in the expression 
\begin{equation} \label{phi(u)}
E=\varphi_2(u_1)\dots\varphi_2(u_m)\varphi_2(u_1^a)\dots\varphi_2(u_m^a). 
\end{equation}
Of course, the element represented by $E$ is in $H_{\sigma^t\omega}$ and 
can be rewritten in the form $b_1^{g_1}\dots b_k^{g_k}$, where $b_i\in 
\{b,c,d\}$ and $g_i\in \{1,a\}$. Let $X_o(E)$ and $X_e(E)$ denote the 
product of the $B$-letters in $E$ preceded by an odd and even number, 
respectively, of $a$'s. Those $B$-letters preceded by odd number of $a$'s 
will appear conjugated by $a$ when we rewrite $E$ in the form 
$b_1^{g_1}\dots b_k^{g_k}$ and those preceded by even number of $a$'s will 
appear without conjugation. It is not difficult to see that we have either 
$X_o(E)=X_e(E)=1$ or $X_o(E)=X_e(E)=b$. Indeed, if the number of $a$'s in 
the expression $\varphi_2(u_1)\dots\varphi_2(u_m)$ is odd, i.e, 
$\varphi_2(u_1)\dots\varphi_2(u_m)$ is not in $H$ then both the number of 
$c$ factors and the number of $d$ factors in both $X_o(E)$ and $X_e(E)$ 
are even so they cancel out and the number of $b$ factors in  $X_o(E)$ and 
$X_e(E)$ is equal, so that their product is 1 or $b$. Similarly, if the 
number of $a$'s in the expression $\varphi_2(u_1)\dots\varphi_2(u_m)$ is 
even then the number of $b$ factors in both $X_o(E)$ and $X_e(E)$ is even 
so the $b$'s cancel out and the number of $c$ factors and $d$ factors in 
$X_o(E)$ and $X_e(E)$ is equal and even so their product is 1 or $b$. 

Considering the extra $b$ in the expression for $h$ we may suppose, up to 
a permutation of the indices $o$ and $e$, that $X_o(h)=b$ and $X_e(h)=1$. 
Then
\begin{align*}
  \chi(\varphi_1(h))&=(1,b),\\
  \chi(\varphi_2(h))&= \begin{cases} (a,1), \quad\text{if }b\not\in K_{t+1} \\
    (1,1), \quad\text{if } b \in K_{t+1} 
  \end{cases}
\end{align*}

The same argument works when we start with $\bar{f}\in
[G_{\sigma^{t+1}\omega},G_{\sigma^{t+1}\omega}]a$ and set
$\bar{h}=\varphi_2(\bar{f}^2)$; we then obtain either $X_o=X_e=1$
whence $\chi(\varphi_i(\bar{h}))=(1,1)$ for $i=1,2$, or
$X_o=X_e=\bar{b}$, where $\bar{b}$ is the only $B$-letter in the
kernel $K_{t+2}$, whence
\[ \chi(\varphi_i(\bar{h})) = 
        \begin{cases} (a,\bar{b}), \quad\text{if } \bar{b} \not\in K_{t+3} \\
                  (1,\bar{b}), \quad\text{if } \bar{b} \in K_{t+3} 
        \end{cases} \quad \text{for } i\in\{1,2\}.
\]

Thus, we enter the graph above at the vertex $(a,b)$, loop $t-1$ times at 
$(a,b)$, then in step $t$ move to the vertex $(1,b)$ and either use two 
consecutive single arrows to exit the graph and lend into an element of 
$G_{\sigma^{t+2}\omega}$ or we proceed to the vertex $(a,1)$, follow 
another double (squaring) arrow to $(1,1)$ and then either leave the graph 
through two consecutive single arrows and lend into an element of 
$G_{\sigma^{t+4}\omega}$ or we lend in an element in 
$G_{\sigma^{t+3}\omega}$ with projection $(a,\bar{b})$. The order of $g$ 
will depend on the exit point, i.e.\ on the number of squarings performed 
and the various length reductions that occurred during the trip through 
the graph. Recall that each time we follow a non-squaring arrow we may 
claim a length reduction by a factor of $1/2$ and in case we follow $r$ 
consecutive squaring arrows (i.e.\ $t=r$) we may claim an additional 
length reduction by a factor of $2/3$. 

As usual, we can now use induction on $n$ and prove the statement 
simultaneously for all $r$-homogeneous $\omega$ (just like in the proof of 
Theorem \ref{r-1/1}).  
\end{proof}

Roughly speaking the previous theorem says that the ratio between the 
number of squarings and the number of halvings performed to calculate the 
order of an element does not exceed $r/2$, i.e.\ in the worst case each 
$r$ squarings are accompanied by at least 2 halvings. Theorem \ref{r-1/1} 
from the previous subsection states, more moderately, that in the worst 
case each $r-1$ raisings to the $p$-th power are accompanied by at least 
one halving step. 

%----------------------------------------------------------------
\subsection{Lower bounds on period growth}
In this subsection we present a construction of words of ``large''
order and ``small'' length, thus providing a lower bound on the period
growth of some regular spinal groups.

The construction in the (proof of the) next theorem generalizes an
unpublished idea of Igor Lysionok, who constructed short words of high
order in the first Grigorchuk group.

Even though the words we construct are far from optimal, they give a 
polynomial lower bound on the degree of period growth in the considered 
cases. More precisely, we show: 

\begin{theorem} \label{thm:periodlb}
  Let $\omega$ be a word in $\Omegah$ and $a\in G_A$ be of order $2$
  and satisfy $a(1)=q$. For all $j\in \N$, set $K_{j,a}=\{b\in
  B|\,\omega_j(b)=a\}$ and
  \[I=\{1\} \cup \{i>1|\,\omega_1=\omega_i\text{ and }\omega_{i-1}\neq\omega_i\}.\]
  Assume that $I$ is infinite, the difference between two consecutive
  indices in $I$ is at most $r$, $K_{1,a} \cap K_{j,a} \neq \emptyset$
  for all $j>1$ and $K_1 \cap K_{j-1,a} \neq \emptyset$ for $j\in
  I,j>1$. Then
  \[ n^{1/(r-1)} \precsim \pi_\omega(n) .\] 
\end{theorem}
\begin{proof}
  Assume that $s+1$ lies in $I$ and that $G_{\sigma^s\omega}$ contains an
  element $g$ satisfying the following conditions:
\begin{enumerate}
  \item $g$ is of order $M$;
  \item $g$ has a representation of length $2k$ with $k$ odd;
    \label{cond:nondiv}
  \item this representation is of the form $ab_1ab_2\dots ab_k$, with
    all the $b_i$ in $K_{s+1,a}$ except for one, which is in $K_{s+1}$. 
  \label{cond:nf}
\end{enumerate}
We shall construct an element $g'$ of $G_\omega$ of order at least $2M$, 
having a representation of length $2(2^{s-1}k+1)$ satisfying Condition 
\ref{cond:nf} (with $1$ instead of $s+1$). 

We restrict our attention to words of the form $a*a*\dots a*$, with $*\in 
B$. A \mathem{word-set} is such a word, but where the $*$'s are non-empty 
subsets of $B$. An \mathem{instance} of a word-set is a word (or group 
element) obtained by choosing an element in each set. 

For all $i\in \N$, there is a map from word-sets in $G_{\sigma^i\omega}$ 
to word-sets in $H_{\sigma^{i-1}\omega}$, defined as follows: 
\begin{gather*}
  \beta_i(b) = b \text{ for all } b \subseteq B;\\
  \beta_i(a) = aK_{i,a}a.
\end{gather*}
Any instance $h'$ of $\beta_i(h)$ satisfies $\psi(h')=(*,\dots,*,h'')$ 
where $h''$ is an instance of $h$. 

Let us now consider $g\in G_{\sigma^s\omega}$ of order $M$. Set first 
$g''=\beta_1\beta_2\dots\beta_s(g)$ and note that 
$\psi_s(g'')=(g,*,\dots,*)$, where each $*$ represents an element in 
$G_{\sigma^s\omega}$ and we are not interested in their actual value. 

Choose $x_1,x_2\in K_{1,a}$ such that $x_1x_2=b_i$, where $b_i$ is the
only letter from $K_1$ in $g''$ and replace $b_i$ by the word-set
$x_1aK_1ax_2$. Note that $K_{i,a}^2=K_i$ so that the choice indicated in 
the previous sentence can be done. Also note that this transformation does 
not change the first coordinate of $\psi_s(g'')$. We claim that $g''$ has 
an instance that, as a word, is a square, say of $g'$. Then since 
$\psi_s((g')^2)=(*,\dots,*,g)$, we will have constructed a word $g'$ 
satisfying the required conditions. 

Let us compute the lengths. Before the substitution of the element in 
$K_1$, $g''$ is a word-set of length $2^{s+1}k$. The substitution of the 
element in $K_1$ increases the length by $4$. Thus, the length of the word 
set $g''$ is $2^{s+1}k+4$. Also, the number of appearances of $K_1$ in 
$g''$ is 1.  

Write $g''$ in 2 lines of length $2^sk+2$. Each line will be of the
form $a*\dots a*$, with the $*$'s elements or subsets of $B$. Our goal
is to choose an instance of $g''$ such that the two lines are
identical, i.e, we want to choose identical elements in each column.
Half of the columns will consist of $a$'s and the other $2^{s-1}k+1$
columns (an odd number of them!)  will have one of the following:
\begin{enumerate}
  \item $K_{1,a}$ and $K_{j,a}$; 
  \item $K_{1,a}$ and $b_i\in K_{s+1,a}=K_{1,a}$ ($k-1$ columns);
  \item $K_{1,a}$, and $x_i$ (two columns); 
  \item $K_1$ and $K_{s,a}$ (one column). 
\end{enumerate}
In each case the two elements in the column can be chosen identical. 

The above construction is a single step in an inductive construction
in which, starting with an element $ab \in G_{\sigma^{t-1}\omega}$
where $t \in I$, $b \in K_1=K_t$, $x$ repetitions of the step give an
element in $G_\omega$ of order at least $2^{x+2}$ and length at most
$2\frac{2^{(r-1)(x+1)}-1}{2^{r-1}-1}$, so that $n^{\frac{1}{r-1}}
\precsim \pi_\omega(n)$.
\end{proof}

Just for the sake of illustration, let us consider an example. Assume
$\omega=\omega_1\omega_2\omega_3\omega_1\dots$ with $\omega_1 \neq
\omega_3$ and start with the element $g=ab_1ab_2ad \in G_{\sigma^3\omega}$ 
where $d \in K_4=K_1$ and $b_1,b_2 \in K_{4,a}=K_{1,a}$. Let $x_1x_2=d$ 
where $x_1,x_2\in K_{1,a}$. We have 
{\allowdisplaybreaks 
\begin{gather*} 
\beta_1\beta_2\beta_3(g)=\beta_1\beta_2(aK_{3,a}ab_1aK_{3,a}ab_2aK_{3,a}ad)=\\ 
=\beta_1(aK_{2,a}aK_{3,a}aK_{2,a}ab_1aK_{2,a}aK_{3,a}aK_{2,a}ab_2aK_{2,a}aK_{3,a}aK_{2,a}ad)=\\ 
=aK_{1,a}aK_{2,a}aK_{1,a}aK_{3,a}aK_{1,a}aK_{2,a}aK_{1,a}ab_1aK_{1,a}aK_{2,a}aK_{1,a}aK_{3,a} 
aK_{1,a}\\* 
aK_{2,a}aK_{1,a}ab_2aK_{1,a}aK_{2,a}aK_{1,a}aK_{3,a}aK_{1,a}aK_{2,a}aK_{1,a}ad. 
\end{gather*}} 
After we replace $d$ by $x_1aK_1ax_2$, write down the word $g''$ in two 
lines and omit the 13 columns consisting of two $a$'s we get 
\[\setlength\arraycolsep{0.3em}\begin{array}{ccccccccccccc}
  K_{1,a} &K_{2,a} &K_{1,a} &K_{3,a} &K_{1,a} &K_{2,a} &K_{1,a} &b_1 
  &K_{1,a} &K_{2,a} &K_{1,a} &K_{3,a} &K_{1,a}\\
  K_{2,a} &K_{1,a} &b_2&K_{1,a} &K_{2,a} &K_{1,a} &K_{3,a} &K_{1,a}
  &K_{2,a} &K_{1,a} &x_1 &K_1&x_2
\end{array}\]
and the entries in the columns are as described in the proof of the
theorem.

If we want to be more specific, assume that we are dealing with a
Grigorchuk $2$-group defined by a sequence that starts with
$\omega=0120\dots$ and $g=abacad$. Then $K_1=\{1,d\}$,
$K_{1,a}=\{b,c\}$, $K_{2,a}=\{b,d\}$, $K_{3,a}=\{c,d\}$, a possible
choice for $x_1$ and $x_2$ is $x_1=b$, $x_2=c$ and an instance of
$g''$ is the square of
\[ g' = ababacacababacababababadac. \] 

It appears the theorem has a lot of assumptions, but all of them 
(except for the existence of $r$) are satisfied, for example, by any 
spinal group with $q=2$. 

Let us point out that the theorem shows that there are regular spinal
groups with at least linear degree of period growth (those defined by
$\omega$ as in the theorem with $r=2$). Moreover, uncountably many
examples can be easily found among the Grigorchuk $2$-groups (any
Grigorchuk $2$-group defined by $0*0*0*0*0*0*\dots$ where the $*$'s
represent arbitrary letters in $\{1,2\}$).

%-----------------------------------------------------------------------------
\section{Final Remarks, Open Questions and Directions}
It is noticeable that most of the conditions (like $r$-homogeneous or
$r$-factorable) put on the defining sequences $\omega$ throughout the
text require appearances of some homomorphisms or subsequences of
homomorphisms in a regular fashion with a frequency that could be
described and bounded uniformly by the number $r$. We might note here
that any initial segment of $\omega$ has no influence on the
asymptotics of the growth functions in spinal groups, so that if the
desired nice behaviour of $\omega$ begins with a little bit of delay
we can still use it. More generally, in some cases we could relax the
``uniformly bounded by $r$'' type of conditions to limit conditions
that describe the density of appearances of the homomorphisms or
subsequences with the desired property.

In the authors' opinion, no spinal group can be finitely presented.
The residual finiteness immediately implies this for the spinal groups
with non-solvable word problem (uncountably many of them). The
solvable word problem case is more involved and interesting and it can
probably be handled in a way similar to the way Grigorchuk $2$-groups
were treated in \cite{grigorchuk:gdegree}.

This would not come as a surprise if one showed that all spinal groups
are branch and just infinite. See~\cite{grigorchuk:jibg} for
definitions and for more constructions of groups of similar flavor.

An interesting direction in the investigation of the growth problems
is introduced by Yuri\u\i\ Leonov in~\cite{leonov:estimation} where he
connects explicitly the word and period growth of some Grigorchuk
$2$-groups. For example, Leonov proves that if $\gamma(n) \precsim
e^{n^\alpha}$, where $0< \alpha <1$, holds for the degree of growth of
the first Grigorchuk group, then $\pi(n) \precsim n^{3\alpha}$ holds
for the degree of period growth. It would be interesting to describe
connections of a similar type in the more general setting of the
present paper.

%-----------------------------------------------------------------------------
\section*{Acknowledgments} 
Both authors are thankful to Pierre de la Harpe and Rostislav Grigorchuk 
for their support, useful suggestions and guidance through the subject. 

In addition, the second author would like to express gratitude to his 
advisor Fernando Guzm\'{a}n for providing infinite listening patience, 
encouragement and plenty of good ideas. Also, thanks to Matthew Brin for a 
partial introduction to the subject. 

%-----------------------------------------------------------------------------
\def\nop#1{}\font\cyr=wncyr8
\providecommand{\bysame}{\leavevmode\hbox to3em{\hrulefill}\thinspace}


\begin{thebibliography}{Bar00b}

\bibitem[Bar98]{bartholdi:upperbd}
Laurent Bartholdi, \emph{The growth of {Grigorchuk}'s torsion group}, Internat.
  Math. Res. Notices \textbf{20} (1998), 1049--1054.

\bibitem[Bar00a]{bartholdi:ggs}
Laurent Bartholdi, \emph{A class of groups acting on rooted trees},
  unpublished, 2000.

\bibitem[Bar00b]{bartholdi:lowerbd}
Laurent Bartholdi, \emph{Lower bounds on the growth of {Grigorchuk}'s torsion
  group}, to appear in Internat. J. Algebra Comput., 2000.

\bibitem[Bas72]{bass:nilpotent}
Hyman Bass, \emph{The degree of polynomial growth of finitely generated
  nilpotent groups}, Proc. London Math. Soc. (3) \textbf{25} (1972), 603--614.

\bibitem[BG99]{bartholdi-g:spectrum}
Laurent Bartholdi and Rostislav~I. Grigorchuk, \emph{On the spectrum of {H}ecke
  type operators related to some fractal groups}, submitted, 1999.

\bibitem[BG00]{bartholdi-g:lie}
Laurent Bartholdi and Rostislav~I. Grigorchuk, \emph{Lie methods in growth of
  groups and groups of finite width}, Computational and Geometric Aspects of
  Modern Algebra (Michael~Atkinson et~al., ed.), London Math. Soc. Lect. Note
  Ser., vol. 275, Cambridge Univ. Press, Cambridge, 2000, pp.~1--27.

\bibitem[Bri98]{brin:igt}
Matthew~G. Brin, \emph{Groups acting on 1-dimensional spaces}, notes for a
  course at SUNY-B, Binghamton; available at
  \texttt{http://math.binghamton.edu/matt/index.html}, 1998.

\bibitem[FG91]{gupta-f:growth2}
Jacek Fabrykowski and Narain~D. Gupta, \emph{On groups with sub-exponential
  growth functions. {II}}, J. Indian Math. Soc. (N.S.) \textbf{56} (1991),
  no.~1-4, 217--228.

\bibitem[GC71]{gecseg-c:ata}
Ferenc G{\'e}cseg and B{\'e}la Cs{\'a}k{\'a}ny, \emph{Algebraic theory of
  automata}, Akademiami Kiado, Budapest, 1971.

\bibitem[GH90]{ghys-h:gromov}
{\'E}tienne Ghys and Pierre~{de la} Harpe, \emph{Sur les groupes hyperboliques
  d'apr\`es {Mikhael Gromov}}, Progress in Mathematics, vol.~83, Birkh{\"a}user
  Boston Inc., Boston, MA, 1990, Papers from the Swiss Seminar on Hyperbolic
  Groups held in Bern, 1988.

\bibitem[Gri80]{grigorchuk:burnside}
Rostislav~I. Grigorchuk, \emph{On {B}urnside's problem on periodic groups},
  Funktsional. Anal. i Prilozhen. \textbf{14} (1980), no.~1, 53--54, English
  translation: {Functional Anal. Appl. \textbf{14} (1980), 41--43}.

\bibitem[Gri83]{grigorchuk:growth}
Rostislav~I. Grigorchuk, \emph{On the {M}ilnor problem of group growth}, Dokl.
  Akad. Nauk SSSR \textbf{271} (1983), no.~1, 30--33.

\bibitem[Gri84]{grigorchuk:gdegree}
Rostislav~I. Grigorchuk, \emph{Degrees of growth of finitely generated groups
  and the theory of invariant means}, Izv. Akad. Nauk SSSR Ser. Mat.
  \textbf{48} (1984), no.~5, 939--985, English translation: {Math. USSR-Izv.
  \textbf{25} (1985), no.~2, 259--300}.

\bibitem[Gri85]{grigorchuk:pgps}
Rostislav~I. Grigorchuk, \emph{Degrees of growth of $p$-groups and torsion-free
  groups}, Mat. Sb. (N.S.) \textbf{126(168)} (1985), no.~2, 194--214, 286.

\bibitem[Gri89]{grigorchuk:hp}
Rostislav~I. Grigorchuk, \emph{On the {H}ilbert-{P}oincar\'e series of graded
  algebras that are associated with groups}, Mat. Sb. \textbf{180} (1989),
  no.~2, 207--225, 304, English translation: {Math. USSR-Sb. \textbf{66}
  (1990), no.~1, 211--229}.

\bibitem[Gri99]{grigorchuk:bath}
Rostislav~I. Grigorchuk, \emph{On the system of defining relations and the
  {S}chur multiplier of periodic groups generated by finite automata}, Groups
  St. Andrews 1997 in Bath, I, Cambridge Univ. Press, Cambridge, 1999,
  pp.~290--317.

\bibitem[Gri00]{grigorchuk:jibg}
Rostislav~I. Grigorchuk, \emph{Just infinite branched groups}, Horizons in
  Profinite Groups (Dan Segal, Markus P.~F. du~Sautoy, and Aner Shalev, eds.),
  Birkha{\"u}ser, Basel, 2000, pp.~121--179.

\bibitem[Gro81]{gromov:nilpotent}
Mikhael Gromov, \emph{Groups of polynomial growth and expanding maps}, Inst.
  Hautes {\'E}tudes Sci. Publ. Math. (1981), no.~53, 53--73.

\bibitem[GS95]{gessel-s:handbook}
Ira~M. Gessel and Richard~P. Stanley, \emph{Handbook of combinatorics},
  ch.~Algebraic Enumeration, pp.~1021--1061, Elsevier, Amsterdam, 1995.

\bibitem[Gui70]{guivarch:poly1}
Yves Guivarc'h, \emph{Groupes de {L}ie \`a croissance polynomiale}, C. R. Acad.
  Sci. Paris S{\'e}r. A-B \textbf{271} (1970), A237--A239.

\bibitem[Har00]{harpe:cgt}
Pierre~{de la} Harpe, \emph{Topics in geometric group theory}, University of
  Chicago Press, 2000.

\bibitem[Leo98]{leonov:pont}
Yuri\u\i~G. Leonov, \emph{On growth function for some torsion residually finite
  groups}, International Conference dedicated to the 90th Anniversary of
  L.S.Pontryagin (Moscow), vol. Algebra, Steklov Mathematical Institute,
  September 1998, pp.~36--38.

\bibitem[Leo99]{leonov:estimation}
Yuri\u\i~G. Leonov, \emph{On precisement of estimation of periods' growth for
  {Grigorchuk}'s $2$-groups}, unpublished, 1999.

\bibitem[Mil68]{milnor:5603}
John~W. Milnor, \emph{Problem 5603}, Amer. Math. Monthly \textbf{75} (1968),
  685--686.

\bibitem[MP99]{muchnik-p:growth}
Roman Muchnik and Igor Pak, \emph{On growth of {Grigorchuk} groups}, preprint,
  1999.

\bibitem[Roz86]{rozhkov:aleshin}
Alexander~V. Rozhkov, \emph{{\cyr K Teorii Grupp Aleshinskogo Tipa} (russian)},
  Mat. Zametki \textbf{40} (1986), no.~5, 572--589.

\bibitem[Tit72]{tits:linear}
Jacques Tits, \emph{Free subgroups in linear groups}, J. Algebra \textbf{20}
  (1972), 250--270.

\end{thebibliography}
\end{document}